\documentclass[11pt]{article}

\usepackage[T1]{fontenc}
\usepackage{amsmath,amssymb,amsthm,mathrsfs}
\usepackage{newtxtext,newtxmath}
\usepackage{dsfont}
\usepackage[margin=1in]{geometry}
\usepackage[colorlinks=true,linkcolor=blue,citecolor=blue,urlcolor=blue]{hyperref}

\allowdisplaybreaks[3]
\newtheorem{theorem}{Theorem}[section]
\newtheorem{proposition}[theorem]{Proposition}
\newtheorem{lemma}[theorem]{Lemma}
\newtheorem{corollary}[theorem]{Corollary}

\newtheorem{remark}[theorem]{Remark}

\numberwithin{equation}{section}

\title{\bfseries\Large Rate of convergence for the Wiener--It\^o chaos in $W_p$ for $1\leq p<\infty$}
\author{Zizhan Chen \and Xiang-Dong Li}
\date{}

\begin{document}

\maketitle

\begin{abstract}
Inspired by previous work of Nualart, Ledoux, Peccati, Nourdin, and others, we establish quantitative fourth moment theorems for random variables in a fixed Wiener--It\^o chaos and for sums of two Wiener--It\^o chaoses in the $W_p$-Wasserstein distance for $1\leq p<\infty$.
\end{abstract}

\noindent\textbf{Keywords.}
Fourth moment theorem; $W_p$-Wasserstein distance; Malliavin calculus; Wiener--It\^o chaos; Stein kernel; coarea formula.

\section{Introduction}

The central limit theorem is a fundamental topic in probability theory and statistics, with a long history of important developments. The Lindeberg--Feller central limit theorem \cite{FELLER71} is a general classical result for sums of independent, non-identically distributed random variables. For Gaussian approximation of nonlinear functionals, however, a general central limit theorem may require convergence of all moments. In the framework of Wiener--It\^o chaos, Nualart and Peccati \cite{NP05} discovered a remarkable simplification. More precisely, let $X=\{X(h):h\in\mathbb{H}\}$ be an isonormal Gaussian process over a real separable Hilbert space $\mathbb{H}$. Fix an integer $q\geq2$, and let
\begin{equation*}
	F_m=I_q(f_m),
	\qquad
	f_m\in\mathbb{H}^{\odot q},
	\qquad
	q!\left\lVert f_m\right\rVert_{\mathbb{H}^{\otimes q}}^2=1,
\end{equation*}
where $\mathbb{H}^{\odot q}$ is the tensor product of $\mathbb{H}$. Denote by $\mu_m$ the distribution of $F_m$ and by $\gamma$ the standard Gaussian distribution.
Then the following assertions are equivalent:
\begin{itemize}
	\item $\mathbb{E}[F_m^4]\to3$ as $m\to\infty$;
	\item $\mu_m$ converges weakly to $\gamma$.
\end{itemize}
This result, now known as the fourth moment theorem, shows that within a fixed Wiener--It\^o chaos the fourth moment alone detects asymptotic normality.

Peccati and Tudor \cite{PT05} subsequently established the corresponding multidimensional principle.  More precisely, let
\begin{equation*}
	F_m=\left(F_{m,1},\ldots,F_{m,d}\right),
	\qquad
	F_{m,i}=I_{q_i}(f_{m,i}),
\end{equation*}
where $q_1,\ldots,q_d$ are fixed.  Denote by $\mu_m$ the distribution of $F_m$.  If the covariance matrices of $F_m$ converge to a matrix $C$, denote by $\gamma_C$ the centered Gaussian distribution with covariance matrix $C$.  Then componentwise convergence of $F_{m,i}$ to the corresponding Gaussian marginals is equivalent to the weak convergence
\begin{equation*}
	\mu_m\Rightarrow\gamma_C.
\end{equation*}
Thus, after the covariances have been identified, no additional mixed-moment assumptions are needed.  Fourth moment theorems for complex multiple Wiener--It\^o integrals were later obtained in related settings; see, for example, Chen and Liu \cite{CL17} and Campese \cite{CAM13}.

The Malliavin--Stein method turns these qualitative statements into quantitative estimates.  Let $F=I_q(f)$ satisfy $\mathbb{E}[F^2]=1$, and denote by $\mu$ the distribution of $F$ and by $\gamma$ the standard Gaussian distribution.  Nourdin and Peccati \cite{NP09,NP12,NP15} proved
\begin{eqnarray*}
	d_{\mathrm{TV}}\left(\mu,\gamma\right)
	&\leq&
	2\sqrt{\frac{q-1}{3q}\left(\mathbb{E}[F^4]-3\right)},\\
	d_{\mathrm K}\left(\mu,\gamma\right)
	&\leq&
	\sqrt{\frac{q-1}{3q}\left(\mathbb{E}[F^4]-3\right)},\\
	W_1\left(\mu,\gamma\right)
	&\leq&
	\sqrt{\frac{2q-2}{3\pi q}\left(\mathbb{E}[F^4]-3\right)}.
\end{eqnarray*}
For probability measures $\mu$ and $\nu$ on $\mathbb{R}^d$, the total variation distance is
\begin{equation*}
	d_{\mathrm{TV}}(\mu,\nu)
	:=
	\sup_{A\in\mathscr B(\mathbb{R}^d)}\left\lvert \mu(A)-\nu(A)\right\rvert.
\end{equation*}
In dimension one, the Kolmogorov distance is
\begin{equation*}
	d_{\mathrm K}(\mu,\nu)
	:=
	\sup_{x\in\mathbb{R}}\left\lvert \mu((-\infty,x])-\nu((-\infty,x])\right\rvert,
\end{equation*}
and, for $p\geq1$, the $W_p$-Wasserstein distance is
\begin{equation*}
	W_p(\mu,\nu)
	:=
	\inf_{\pi\in\Pi(\mu,\nu)}
	\left(
	\int_{\mathbb{R}^d\times\mathbb{R}^d}\left\lvert x-y\right\rvert^p\,d\pi(x,y)
	\right)^{\frac{1}{p}}.
\end{equation*}
Here $\Pi(\mu,\nu)$ denotes the set of all couplings of $\mu$ and $\nu$.  The optimal asymptotic rate in total variation was identified in \cite{NP15}.

The multidimensional counterpart was developed by Nourdin, Peccati and R\'eveillac \cite{NPR10}.  Let
\begin{equation*}
	F=(F_1,\ldots,F_d),
	\qquad
	F_i=I_{q_i}(f_i),
\end{equation*}
and let $C$ be a positive definite covariance matrix.  Denote by $\mu$ the distribution of $F$ and by $\gamma_C$ the centered Gaussian measure with covariance matrix $C$.  Their estimate is expressed through the Malliavin matrix
\begin{equation}\label{eq:Gamma-matrix-intro}
	\Gamma_{ij}(F)
	:=
	\left\langle DF_j,-DL^{-1}F_i\right\rangle_{\mathbb{H}}.
\end{equation}
In particular,
\begin{equation*}
	W_1\left(\mu,\gamma_C\right)
	\leq
	\left\lVert C^{\frac{1}{2}}\right\rVert_{\mathrm{op}}
	\left\lVert C^{-1}\right\rVert_{\mathrm{op}}
	\left(
	\sum_{i,j=1}^d
	\mathbb{E}\left[\left\lvert \Gamma_{ij}(F)-C_{ij}\right\rvert^2\right]
	\right)^{\frac{1}{2}}.
\end{equation*}
Further refinements, including optimal multidimensional rates, were obtained by Campese \cite{CAM13}.

Fourth moment phenomena have also been investigated outside a single homogeneous Gaussian chaos.  Poisson-space analogues were established, under increasingly general assumptions, by D\"obler, Vidotto and Zheng \cite{DVZ18}.  In free probability, Kemp, Nourdin, Peccati and Speicher \cite{KNPS12} proved the free fourth moment theorem for Wigner chaos and obtained an explicit quantitative estimate for finite Wigner-chaos expansions.  Subsequent quantitative extensions were developed by Bourguin and Campese \cite{BC17} and C\'ebron \cite{CEB21}.

More recently, attention has shifted to sums of Wiener--It\^o chaoses. Basse-O'Connor, Kramer-Bang and Svendsen \cite{BKS25} proved a fourth moment theorem for a sum of two multiple Wiener--It\^o integrals whose orders have opposite parity. Under the unit-variance normalization, let $\mu$ denote the distribution of $I_q(f)+I_r(g)$. With the help of Chat-GPT-5, Diez, da Maia and Nourdin \cite{DMN25} obtained the explicit estimate
\begin{equation*}
	d_{\mathrm{TV}}\left(\mu,\gamma\right)
	\leq
	\sqrt{6}\sqrt{\kappa_4\left(I_q(f)+I_r(g)\right)}.
\end{equation*}
For every centered random variable $Y\in L^4$, we write
\begin{equation*}
	\kappa_4(Y):=\mathbb{E}[Y^4]-3\mathbb{E}[Y^2]^2.
\end{equation*}

Motivated by these results, the first purpose of this paper is to obtain quantitative fourth moment bounds in the $W_p$-Wasserstein distance for every $1\leq p<\infty$. For $p>1$, convergence in $W_p$ is stronger than convergence in $W_1$. The second purpose is to give an alternative normal-approximation argument, without Stein's method, for an orthogonally diagonalizable Wiener--It\^o chaos of arbitrary fixed order. This argument combines Hermite small-ball estimates, negative moments of the Malliavin covariance and the coarea formula.

We now state the main results.  Define
\begin{equation}\label{eq:Cp}
	C_p
	:=
	\left(\int_{\mathbb{R}}\left\lvert x\right\rvert^p\,d\gamma(x)\right)^{\frac{1}{p}}
	=
	\sqrt{2}
	\left(
	\frac{\Gamma\left(\frac{p+1}{2}\right)}{\sqrt\pi}
	\right)^{\frac{1}{p}},
\end{equation}
and, for chaos orders $q,r\geq1$, define
\begin{equation}\label{eq:B-constant}
	B_{p,q,r}
	:=
	\begin{cases}
	1,
	&1\leq p\leq2,\\[1mm]
	\left(\displaystyle\sum_{k=0}^{q+r-2}(p-1)^k\right)^{\frac{1}{2}},
	&p>2,
	\end{cases}
\end{equation}
and
\begin{equation}\label{eq:C-constant}
	C_{p,q,r}
	:=
	\begin{cases}
	1,
	&1\leq p\leq2,\\[1mm]
	\left(\displaystyle\sum_{k=1}^{q+r-2}(p-1)^k\right)^{\frac{1}{2}},
	&p>2.
	\end{cases}
\end{equation}
For the multidimensional estimates, put
\begin{equation}\label{eq:beta}
	\beta(p)
	:=
	\begin{cases}
	\frac{1}{p},&1\leq p<2,\\[1mm]
	2\left(1-\frac{2}{p}\right),&p\geq2.
	\end{cases}
\end{equation}

The first theorem extends the quantitative fourth moment theorem from $W_1$ to the $W_p$-Wasserstein distance for every $1\leq p<\infty$.

\begin{theorem}\label{thm:single-chaos-Wp}
	Let $F=I_q(f)$ with $f\in\mathbb{H}^{\odot q}$.  Denote by $\mu$ the distribution of $F$ and by $\gamma_\sigma$ the centered Gaussian distribution with variance $\sigma^2$, where $\sigma>0$.  Then, for every $p\geq1$,
	\begin{equation}\label{eq:single-first}
		W_p\left(\mu,\gamma_\sigma\right)
		\leq
		\frac{C_pB_{p,q,q}}{\sigma}
		\left(
		\mathbb{E}\left[
		\left(
		\frac{1}{q}\left\lVert DF\right\rVert_{\mathbb{H}}^2-\sigma^2
		\right)^2
		\right]
		\right)^{\frac{1}{2}}.
	\end{equation}
	Moreover, we have
	\begin{equation}\label{eq:single-fourth}
		W_p\left(\mu,\gamma_\sigma\right)
		\leq
		\frac{C_p}{\sigma}
		\left(
		C_{p,q,q}\sqrt{\frac{q-1}{3q}}\sqrt{\kappa_4(F)}
		+
		\left\lvert \mathbb{E}[F^2]-\sigma^2\right\rvert
		\right),
	\end{equation}
	where the constants $B_{p,q,q}, C_{p}$ and $C_{p,q,q}$ are defined in \eqref{eq:Cp}--\eqref{eq:C-constant}.
\end{theorem}
In the cases of $p=1$, we concide the result in \cite{NP09}. 

To state the result for two chaos components, let $q, r$ be two positive integers with opposite parity, set
\begin{equation*}
	F=I_q(f),
	\qquad
	G=I_r(g),
	\qquad
	Z=F+G.
\end{equation*}

\begin{theorem}\label{thm:two-chaos-Wp}
	Let $q\neq r$ be two positive integers with opposite parity and $m=q\vee r$.  Denote by $\mu$ the distribution of $Z$ and by $\gamma_\sigma$ the centered Gaussian distribution with variance $\sigma^2$, where $\sigma>0$.  For every $p\geq1$,
	\begin{eqnarray}
	W_p\left(\mu,\gamma_\sigma\right)
	&\leq&
	\frac{C_pB_{p,m,m}}{\sigma}
	\left(
	\mathbb{E}\left[\left(\Gamma_Z-\sigma^2\right)^2\right]
	\right)^{\frac{1}{2}},
	\label{eq:two-first}\\
	W_p\left(\mu,\gamma_\sigma\right)
	&\leq&
	\frac{C_p}{\sigma}
	\left[
	C_{p,m,m}\sqrt{\frac{3}{2}\kappa_4(Z)}
	+
	\left\lvert \mathbb{E}[Z^2]-\sigma^2\right\rvert
	\right],
	\label{eq:two-fourth}
	\end{eqnarray}
	where $\Gamma_Z=\langle DZ,-DL^{-1}Z\rangle_{\mathbb{H}}$ and the constants $B_{p,q,q}, C_{p}$ and $C_{p,q,q}$ are defined in \eqref{eq:Cp}--\eqref{eq:C-constant}.
	
\end{theorem}

The next two theorems give the corresponding multidimensional estimates.  For a random $d\times d$ matrix $A=(A_{ij})$, write
\begin{equation*}
	\left\lVert A\right\rVert_p
	:=
	\left(\sum_{i,j=1}^d\mathbb{E}\left\lvert A_{ij}\right\rvert^p\right)^{\frac{1}{p}}.
\end{equation*}

\begin{theorem}\label{thm:vector-chaos}
	Let
	\begin{equation*}
		F=(F_1,\ldots,F_d),
		\qquad
		F_i=I_{q_i}(f_i).
	\end{equation*}
	Let $C$ be a positive definite $d\times d$ covariance matrix.  Denote by $\mu$ the distribution of $F$ and by $\gamma_C$ the centered Gaussian measure with covariance matrix $C$.  Then, for every $p\geq1$,
	\begin{equation}\label{eq:vector-exact-first}
	W_p\left(\mu,\gamma_C\right)
	\leq
	C_p d^{\beta(p)}
	\left\lVert C^{\frac{1}{2}}\right\rVert_{\mathrm{op}}
	\left\lVert C^{-1}\right\rVert_{\mathrm{op}}
	\left(
	\sum_{i,j=1}^d
	B_{p,q_i,q_j}^{p}
	\mathbb{E}\left[\left\lvert \Gamma_{ij}(F)-C_{ij}\right\rvert^2\right]^{\frac{p}{2}}
	\right)^{\frac{1}{p}},
	\end{equation}
	where $\Gamma_Z=\langle DZ,-DL^{-1}Z\rangle_{\mathbb{H}}$. 
	Moreover, there exists a constant $D_{p,d,q_1,\ldots,q_d}>0$ such that
	\begin{equation}\label{eq:vector-fourth-intro}
	W_p\left(\mu,\gamma_C\right)
	\leq
	D_{p,d,q_1,\ldots,q_d}
	\left\lVert C^{\frac{1}{2}}\right\rVert_{\mathrm{op}}
	\left\lVert C^{-1}\right\rVert_{\mathrm{op}}
	\left[
	\sum_{i=1}^d\sqrt{\kappa_4(F_i)}
	+
	\sum_{i\neq j}\mathbb{E}[F_i^2]^{\frac{1}{2}}\kappa_4(F_j)^{\frac{1}{4}}
	+
	\sum_{i,j=1}^d\left\lvert \mathbb{E}[F_iF_j]-C_{ij}\right\rvert
	\right],
	\end{equation}
	where the constants $B_{p,q,q}, C_{p}$ and $C_{p,q,q}$ are defined in \eqref{eq:Cp}--\eqref{eq:C-constant}.
\end{theorem}

For $1\leq i\leq d$, let
\begin{equation*}
	Z_i=F_i+G_i,
	\qquad
	F_i=I_{q_i}(f_i),
	\qquad
	G_i=I_{r_i}(g_i),
	\qquad
	q_i\neq r_i,
\end{equation*}
and put
\begin{equation*}
	m_i:=q_i\vee r_i,
	\qquad
	v_i:=\mathbb{E}[Z_i^2].
\end{equation*}
Write $\boldsymbol q=(q_1,\ldots,q_d)$ and $\boldsymbol r=(r_1,\ldots,r_d)$.

\begin{theorem}\label{thm:vector-two-chaos}
	Let $Z=(Z_1,\ldots,Z_d)$ and let $C$ be a positive definite covariance matrix. Denote by $\mu$ the distribution of $Z$ and by $\gamma_C$ the centered Gaussian measure with covariance matrix $C$. Assume that $q_i$ and $r_i$ have opposite parity for every $i$. Then, for every $p\geq1$,
	\begin{equation}\label{eq:vector-two-exact}
	W_p\left(\mu,\gamma_C\right)
	\leq
	C_p d^{\beta(p)}
	\left\lVert C^{\frac{1}{2}}\right\rVert_{\mathrm{op}}
	\left\lVert C^{-1}\right\rVert_{\mathrm{op}}
	\left(
	\sum_{i,j=1}^d
	B_{p,m_i,m_j}^{p}
	\mathbb{E}\left[\left\lvert \Gamma_{ij}(Z)-C_{ij}\right\rvert^2\right]^{\frac{p}{2}}
	\right)^{\frac{1}{p}},
	\end{equation}
	where $\Gamma_Z=\langle DZ,-DL^{-1}Z\rangle_{\mathbb{H}}$ and the constants $B_{p,q,q}, C_{p}$ and $C_{p,q,q}$ are defined in \eqref{eq:Cp}--\eqref{eq:C-constant}.
	Moreover, there exists a constant $D_{p,d,\boldsymbol q,\boldsymbol r}>0$ such that
	\begin{equation}	\label{eq:vector-two-fourth}
	W_p\left(\mu,\gamma_C\right)
	\leq
	D_{p,d,\boldsymbol q,\boldsymbol r}
	\left\lVert C^{\frac{1}{2}}\right\rVert_{\mathrm{op}}
	\left\lVert C^{-1}\right\rVert_{\mathrm{op}}
	\left[
	\sum_{i=1}^d\kappa_4(Z_i)^{\frac{1}{2}}
	+
	\sum_{i\neq j}v_i^{\frac{1}{2}}\kappa_4(Z_j)^{\frac{1}{4}}
	+
	\sum_{i,j=1}^d\left\lvert \mathbb{E}[Z_iZ_j]-C_{ij}\right\rvert
	\right].
	\end{equation}
	
\end{theorem}

We next give an estimate for finite sums of Wiener--It\^o chaoses in the $W_p$-Wasserstein distance.

\begin{theorem}\label{thm:finite-chaos-sum}
	Let $F\in\mathbb D^{1,2}$ be centered with $\mathbb{E}[F^2]=1$.  Denote by $\mu$ the distribution of $F$ and by $\gamma$ the standard Gaussian distribution.  Define
	\begin{equation*}
		\Gamma_F:=\left\langle DF,-DL^{-1}F\right\rangle_{\mathbb{H}}.
	\end{equation*}
	Also set
	\begin{equation*}
		d_2\left(\mu,\gamma\right)
		:=
		\sup_{\substack{h\in C^2(\mathbb{R}),\ h(0)=h'(0)=0\\
		\left\lVert h''\right\rVert_\infty\leq1}}
		\left\lvert 
		\int_\mathbb{R} h\,d\mu
		-
		\int_\mathbb{R} h\,d\gamma
		\right\rvert.
	\end{equation*}
	Then
	\begin{equation}\label{eq:smart-path-defect}
		d_2\left(\mu,\gamma\right)
		\leq
		\frac{1}{2}\mathbb{E}\left\lvert 1-\Gamma_F\right\rvert.
	\end{equation}
	Suppose in addition that
	\begin{equation}\label{eq:finite-chaos-sum}
		F=\sum_{q=1}^{M}I_q(f_q).
	\end{equation}
	For $1\leq\ell\leq2M-2$, set
	\begin{equation}\label{eq:h-ell-finite-sum}
		h_\ell
		:=
		\sum_{\substack{1\leq q,r\leq M,\ 1\leq k\leq q\wedge r\\q+r-2k=\ell}}
		q(k-1)!
		\binom{q-1}{k-1}\binom{r-1}{k-1}
		\left(f_q\widetilde{\otimes}_k f_r\right).
	\end{equation}
	Then
	\begin{equation}\label{eq:finite-sum-exact-variance}
		\mathbb{E}\left[\left(\Gamma_F-1\right)^2\right]
		=
		\sum_{\ell=1}^{2M-2}\ell!\left\lVert h_\ell\right\rVert_{\mathbb{H}^{\otimes\ell}}^2,
	\end{equation}
	and, for every $p\geq1$,
	\begin{equation}\label{eq:finite-sum-Wp}
		W_p\left(\mu,\gamma\right)
		\leq
		C_pC_{p,M,M}
		\left(
		\sum_{\ell=1}^{2M-2}\ell!\left\lVert h_\ell\right\rVert_{\mathbb{H}^{\otimes\ell}}^2
		\right)^{\frac{1}{2}}.
	\end{equation}
\end{theorem}

The final results concern diagonal kernels.  Fix $n\geq2$ and suppose
\begin{equation*}
	f=\sum_{k\geq1}\lambda_ke_k^{\otimes n},
\end{equation*}
where $(e_k)_{k\geq1}$ is orthonormal.  This representation is automatic when $n=2$, by the spectral theorem for compact self-adjoint Hilbert--Schmidt operators.  For $n\geq3$, it is a genuine structural assumption.  In this case
\begin{equation*}
	F=I_n(f)=\sum_{k\geq1}\lambda_kH_n\left(X(e_k)\right),
\end{equation*}
so $F$ is represented as a sum of independent random variables.

\begin{theorem}\label{thm:uniform-inverse}
	Let
	\begin{equation*}
		F_m
		=
		\sum_{k\geq1}\lambda_{m,k}H_n(\xi_{m,k}),
		\qquad
		n!\sum_{k\geq1}\lambda_{m,k}^2=1,
	\end{equation*}
	where $(\xi_{m,k})_{k\geq1}$ is an independent standard Gaussian family for every $m$.  If $\kappa_4(F_m)\to0$, then, for every $q>0$, there exists $m_0=m_0(q)$ such that
	\begin{equation}\label{eq:uniform-negative}
		\sup_{m\geq m_0}\mathbb{E}\left[\gamma_{F_m}^{-q}\right]<\infty,
	\end{equation}
	where $\gamma_{F_m}=\left\lVert DF_m\right\rVert_{\mathbb{H}}^2$.  Moreover,
	\begin{equation}\label{eq:inverse-limit}
		\mathbb{E}\left[\gamma_{F_m}^{-q}\right]
		\longrightarrow
		n^{-q}.
	\end{equation}
\end{theorem}

\begin{theorem}\label{thm:diagonal-non-Stein}
	Let $F$ be a normalized diagonal $n$th-chaos random variable.  Denote by $\mu_F$ the distribution of $F$ and by $\gamma$ the standard Gaussian distribution.  Assume that $F$ is smoothly nondegenerate.  For every $p\geq1$, $\eta>0$, and $q>1$,
	\begin{equation}	\label{eq:diagonal-main}
	W_p\left(\mu_F,\gamma\right)
	\leq
	C_{p,\eta,n,q}
	\left(
	1+
	\left(\mathbb{E}\left[\gamma_F^{-q}\right]\right)^{\frac{1}{q}}
	\right)^{
	\frac{\eta}{2p(p+\eta+\frac{1}{2})}}
	\kappa_4(F)^{
	\frac{\eta}{4p(p+\eta+\frac{1}{2})}}.
	\end{equation}
\end{theorem}

However, this convergence rate is not as sharp as the result obtained in \cite{NP09} and in Theorem \ref{thm:single-chaos-Wp}.
\begin{corollary}\label{cor:diagonal-convergence}
	Under the assumptions of Theorem~\ref{thm:uniform-inverse}, assume in addition that every $F_m$ has infinitely many nonzero diagonal coefficients.  Denote by $\mu_m$ the distribution of $F_m$ and by $\gamma$ the standard Gaussian distribution.  Then, for every $p\geq1$ and $\eta>0$, there exist $m_0$ and $C_{p,\eta,n}>0$ such that
	\begin{equation}\label{eq:diagonal-sequence-rate}
		W_p\left(\mu_m,\gamma\right)
		\leq
		C_{p,\eta,n}
		\kappa_4(F_m)^{
		\frac{\eta}{4p(p+\eta+\frac{1}{2})}},
		\qquad m\geq m_0.
	\end{equation}
	In particular, $\kappa_4(F_m)\to0$ implies $W_p\left(\mu_m,\gamma\right)\to0$ for every finite $p$.
\end{corollary}

The rest of the paper is organized as follows.  Section~2 recalls the necessary material on Wiener--It\^o chaos, Malliavin operators, hypercontractivity, Wasserstein distance, and Stein kernels.  Section~3 proves the finite-chaos-sum estimate and the scalar and vector-valued $W_p$ bounds.  Section~4 treats orthogonally diagonalizable chaoses of arbitrary fixed order, establishes the inverse-moment estimate, and proves the non-Stein normal-approximation theorem by means of the coarea formula and the classical Berry--Esseen theorem.

\section{Preliminaries}

\subsection{Wiener--It\^o chaos and Malliavin calculus}

Throughout the paper,
\begin{equation*}
	\mathbb{H}=L^2\left(\mathbb{R}_+\right),
	\qquad
	X(h)=\int_0^\infty h(t)\,dW_t,
\end{equation*}
where $W$ is a standard Brownian motion.  For $q\geq1$ and $f\in\mathbb{H}^{\odot q}$, let $I_q(f)$ denote the $q$th multiple Wiener--It\^o integral.  If $q,r\geq1$, $f\in\mathbb{H}^{\odot q}$, and $g\in\mathbb{H}^{\odot r}$, then the isometry and orthogonality relations are
\begin{equation}\label{eq:isometry}
	\mathbb{E}\left[I_q(f)I_r(g)\right]
	=
	\begin{cases}
	q!\left\langle f,g\right\rangle_{\mathbb{H}^{\otimes q}},&q=r,\\
	0,&q\neq r.
	\end{cases}
\end{equation}

For $f\in\mathbb{H}^{\odot q}$, $g\in\mathbb{H}^{\odot r}$, and $0\leq k\leq q\wedge r$, the $k$th contraction is denoted by $f\otimes_k g$, and its symmetrization by $f\widetilde{\otimes}_k g$.  The product formula is
\begin{equation}\label{eq:product-formula}
	I_q(f)I_r(g)
	=
	\sum_{k=0}^{q\wedge r}
	k!\binom{q}{k}\binom{r}{k}
	I_{q+r-2k}\left(f\widetilde{\otimes}_k g\right).
\end{equation}

Every $F\in L^2(\Omega)$ has a unique chaos expansion
\begin{equation}\label{eq:chaos-expansion}
	F=\mathbb{E}[F]+\sum_{q=1}^{\infty}I_q(f_q).
\end{equation}
On its natural domain, the Malliavin derivative and the Ornstein--Uhlenbeck generator are given by
\begin{eqnarray}
	D_tF
	&=&
	\sum_{q=1}^{\infty}qI_{q-1}\left(f_q(\cdot,t)\right),
	\label{eq:D-chaos}\\
	LF
	&=&
	-\sum_{q=1}^{\infty}qI_q(f_q),
	\qquad
	L^{-1}F=-\sum_{q=1}^{\infty}\frac{1}{q}I_q(f_q)
	\label{eq:L-chaos}
\end{eqnarray}
when $\mathbb{E}[F]=0$.  In particular, the integral by parts formula holds
\begin{equation}\label{eq:ibp-basic}
	\mathbb{E}[FG]
	=
	\mathbb{E}\left[\left\langle DG,-DL^{-1}F\right\rangle_{\mathbb{H}}\right]
\end{equation}
for centered $F\in L^2$ and suitable $G\in\mathbb D^{1,2}$.

Throughout Sections~2 and~3, $p$ is reserved for the Wasserstein or $L^p$ exponent, $q$ and $r$ denote chaos orders, and $k$ denotes a contraction index.  In Section~\ref{sec:diagonal-chaos}, $n$ denotes the fixed order of the diagonal chaos, and the auxiliary integrability exponents are specified locally.

Nelson's hypercontractivity inequality states that, for $I_q(f)$ and $p\geq2$,
\begin{equation}\label{eq:hypercontractivity}
	\left\lVert I_q(f)\right\rVert_{L^p(\Omega)}
	\leq
	(p-1)^{\frac{q}{2}}\left\lVert I_q(f)\right\rVert_{L^2(\Omega)}.
\end{equation}

The following lemma records the precise finite-chaos comparison that will be used repeatedly.  The constants $B_{p,q,r}$ and $C_{p,q,r}$ are defined in \eqref{eq:B-constant} and \eqref{eq:C-constant}.

\begin{lemma}\label{lem:finite-hypercontractivity}
	Let
	\begin{equation*}
		U=\sum_{\ell=0}^{q+r-2}I_\ell(h_\ell),
	\end{equation*}
	where $h_\ell\in\mathbb{H}^{\odot\ell}$.  Then, for every $p\geq1$,
	\begin{equation}\label{eq:finite-hypercontractivity}
		\left\lVert U\right\rVert_{L^p}
		\leq
		B_{p,q,r}\left\lVert U\right\rVert_{L^2}.
	\end{equation}
	If $\mathbb{E}[U]=0$, then
	\begin{equation}\label{eq:centered-finite-hypercontractivity}
		\left\lVert U\right\rVert_{L^p}
		\leq
		C_{p,q,r}\left\lVert U\right\rVert_{L^2}.
	\end{equation}
	Consequently, if $F=I_q(f)$, $G=I_r(g)$, and $c\in\mathbb{R}$, then
	\begin{equation}\label{eq:Malliavin-hypercontractivity}
		\left\Vert
		\left\langle DF,-DL^{-1}G\right\rangle_{\mathbb{H}}-c
		\right\Vert_{L^p}
		\leq
		B_{p,q,r}
		\left\Vert
		\left\langle DF,-DL^{-1}G\right\rangle_{\mathbb{H}}-c
		\right\Vert_{L^2}.
	\end{equation}
	For $c=\mathbb{E}[FG]$, the constant $B_{p,q,r}$ in \eqref{eq:Malliavin-hypercontractivity} may be replaced by $C_{p,q,r}$.
\end{lemma}

\begin{proof}
	If $1\leq p\leq2$, the monotonicity of $L^p$-norms on a probability space gives
	\begin{equation*}
		\left\lVert U\right\rVert_{L^p}\leq\left\lVert U\right\rVert_{L^2}.
	\end{equation*}
	This proves both assertions in this range because the two constants equal one.  Suppose now that $p>2$.  By Minkowski's inequality and \eqref{eq:hypercontractivity},
	\begin{equation*}
		\left\lVert U\right\rVert_{L^p}
		\leq
		\sum_{\ell=0}^{q+r-2}(p-1)^{\frac{\ell}{2}}
		\left\lVert I_\ell(h_\ell)\right\rVert_{L^2}.
	\end{equation*}
	Cauchy--Schwarz in the finite sum yields
	\begin{eqnarray*}
		\left\lVert U\right\rVert_{L^p}\leq 
		\left(\sum_{\ell=0}^{q+r-2}(p-1)^\ell\right)^{\frac{1}{2}}
		\left(\sum_{\ell=0}^{q+r-2}
		\left\lVert I_\ell(h_\ell)\right\rVert_{L^2}^2\right)^{\frac{1}{2}}=
		B_{p,q,r}\left\lVert U\right\rVert_{L^2},
	\end{eqnarray*}
	where the last identity follows from orthogonality of distinct Wiener--It\^o chaoses.  If $\mathbb{E}[U]=0$, then $I_0(h_0)=0$, so the same calculation begins at $\ell=1$ and gives \eqref{eq:centered-finite-hypercontractivity}.

	Finally, the product formula shows that
	\begin{equation*}
		\left\langle DF,-DL^{-1}G\right\rangle_{\mathbb{H}}-c
		\in
		\bigoplus_{\ell=0}^{q+r-2}\mathcal H_\ell.
	\end{equation*}
	Thus \eqref{eq:Malliavin-hypercontractivity} follows from \eqref{eq:finite-hypercontractivity}.  The integration-by-parts identity \eqref{eq:ibp-basic} gives
	\begin{equation*}
		\mathbb{E}\left[\left\langle DF,-DL^{-1}G\right\rangle_{\mathbb{H}}\right]
		=
		\mathbb{E}[FG].
	\end{equation*}
	Hence the random variable is centered when $c=\mathbb{E}[FG]$, and \eqref{eq:centered-finite-hypercontractivity} applies.
\end{proof}

Unless a constant is explicitly defined, $C$ denotes a positive universal constant, while the subscripts of a generic constant record the parameters on which it may depend; its value may change from line to line.  The quantities $C_p$, $B_{p,q,r}$, and $C_{p,q,r}$ defined in \eqref{eq:Cp}, \eqref{eq:B-constant}, and \eqref{eq:C-constant}, respectively, retain their stated values.

\subsection{Wasserstein distance and Stein kernels}

For probability measures $\mu$ and $\nu$ on $\mathbb{R}^d$ with finite $p$th moments,
\begin{equation*}
	W_p(\mu,\nu)
	:=
	\inf_{\pi\in\Pi(\mu,\nu)}
	\left(
	\int_{\mathbb{R}^d\times\mathbb{R}^d}\left\lvert x-y\right\rvert^p\,d\pi(x,y)
	\right)^{\frac{1}{p}}.
\end{equation*}

Let $\nu$ be a centered probability measure on $\mathbb{R}^d$. Ledoux, Nourdin and Peccati \cite{LNP15} introduced the Stein kernel for the probability measure.  A matrix-valued function $\tau_\nu$ is a Stein kernel for $\nu$ if
\begin{equation}\label{eq:vector-Stein-kernel}
	\int_{\mathbb{R}^d}x\,\varphi(x)\,d\nu(x)
	=
	\int_{\mathbb{R}^d}\tau_\nu(x)\nabla\varphi(x)\,d\nu(x)
\end{equation}
for every smooth scalar function $\varphi$ for which the two sides are finite.  Let $\gamma_d=\mathcal{N}_d(0,I_d)$.  With $C_p$ as in \eqref{eq:Cp}, Proposition 3.4 of Ledoux, Nourdin and Peccati \cite{LNP15} gives
\begin{equation}\label{eq:LNP-Wp}
	W_p(\nu,\gamma_d)
	\leq
	C_p d^{\alpha(p)}
	\left(
	\sum_{i,j=1}^d
	\int_{\mathbb{R}^d}\left\lvert \tau_\nu^{ij}(x)-\delta_{ij}\right\rvert^p\,d\nu(x)
	\right)^{\frac{1}{p}},
\end{equation}
where
\begin{equation}\label{eq:alpha}
	\alpha(p)
	=
	\begin{cases}
	1-\frac{1}{p},&1\leq p<2,\\[1mm]
	1-\frac{2}{p},&p\geq2.
	\end{cases}
\end{equation}

For a centered scalar Wiener functional $F\in\mathbb D^{1,2}$, define
\begin{equation}\label{eq:GammaF}
	\Gamma_F:=\left\langle DF,-DL^{-1}F\right\rangle_{\mathbb{H}}.
\end{equation}
By \eqref{eq:ibp-basic},
\begin{equation}\label{eq:scalar-Stein-kernel}
	\tau_F(x):=\mathbb{E}\left[\Gamma_F\mid F=x\right]
\end{equation}
is a Stein kernel for the law of $F$, and
\begin{equation*}
	\mathbb{E}[\Gamma_F]=\mathbb{E}[F^2].
\end{equation*}

\section{Stein-kernel bounds and finite chaos sums}\label{sec:finite-sums}
\subsection{Proof of Theorem~\ref{thm:single-chaos-Wp} and \ref{thm:two-chaos-Wp}}

We first recall the standard fourth-cumulant estimate.  If $F=I_q(f)$, then Nourdin and Peccati \cite{NP12} proved the following useful inequality
\begin{equation}\label{eq:key-fixed-chaos}
	\mathbb{E}\left[
	\left(
	\frac{1}{q}\left\lVert DF\right\rVert_{\mathbb{H}}^2-\mathbb{E}[F^2]
	\right)^2
	\right]
	\leq
	\frac{q-1}{3q}\,\kappa_4(F).
\end{equation}

Now we give the proof of Theorem~\ref{thm:single-chaos-Wp}.

\begin{proof}
	Set $Y=\dfrac{F}{\sigma}$ and denote by $\mu_\sigma$ the distribution of $Y$.  Since the map $x\mapsto\sigma x$ sends $\mu_\sigma$ to $\mu$ and sends $\gamma$ to $\gamma_\sigma$, the scaling property of the Wasserstein distance gives
	\begin{equation*}
		W_p\left(\mu,\gamma_\sigma\right)
		=
		\sigma W_p\left(\mu_\sigma,\gamma\right).
	\end{equation*}
	Put
	\begin{equation*}
		\Gamma_F
		=
		\left\langle DF,-DL^{-1}F\right\rangle_{\mathbb{H}}
		=
		\frac{1}{q}\left\lVert DF\right\rVert_{\mathbb{H}}^2.
	\end{equation*}
	By \eqref{eq:scalar-Stein-kernel}, a Stein kernel for $Y$ is
	\begin{equation*}
		\tau_Y(y)
		=
		\frac{1}{\sigma^2}
		\mathbb{E}\left[\Gamma_F\mid F=\sigma y\right].
	\end{equation*}
	The one-dimensional case of \eqref{eq:LNP-Wp} and conditional Jensen's inequality therefore yield
	\begin{eqnarray*}
		W_p\left(\mu,\gamma_\sigma\right)
		&=&
		\sigma W_p\left(\mu_\sigma,\gamma\right)\\
		&\leq&
		\sigma C_p
		\left(
		\mathbb{E}\left[\left\lvert \tau_Y(Y)-1\right\rvert^p\right]
		\right)^{\frac{1}{p}}\\
		&\leq&
		\frac{C_p}{\sigma}
		\left(
		\mathbb{E}\left[\left\lvert \Gamma_F-\sigma^2\right\rvert^p\right]
		\right)^{\frac{1}{p}}.
	\end{eqnarray*}
	The product formula shows that
	\begin{equation*}
		\Gamma_F-\sigma^2
		\in
		\bigoplus_{\ell=0}^{2q-2}\mathcal H_\ell.
	\end{equation*}
	Consequently, Lemma~\ref{lem:finite-hypercontractivity} gives
	\begin{equation*}
		\left(
		\mathbb{E}\left[\left\lvert \Gamma_F-\sigma^2\right\rvert^p\right]
		\right)^{\frac{1}{p}}
		\leq
		B_{p,q,q}
		\left(
		\mathbb{E}\left[\left(\Gamma_F-\sigma^2\right)^2\right]
		\right)^{\frac{1}{2}},
	\end{equation*}
	which proves \eqref{eq:single-first}.

	For the second estimate, use $\mathbb{E}[\Gamma_F]=\mathbb{E}[F^2]$ and write
	\begin{equation*}
		\Gamma_F-\sigma^2
		=
		\left(\Gamma_F-\mathbb{E}[F^2]\right)
		+
		\left(\mathbb{E}[F^2]-\sigma^2\right),
	\end{equation*}
	where $\Gamma_F-\mathbb{E}[F^2]$ is centered.  Minkowski's inequality and \eqref{eq:centered-finite-hypercontractivity} imply
	\begin{equation*}
		\left(
		\mathbb{E}\left[\left\lvert \Gamma_F-\sigma^2\right\rvert^p\right]
		\right)^{\frac{1}{p}}
		\leq
		C_{p,q,q}
		\left(
		\mathbb{E}\left[\left(\Gamma_F-\mathbb{E}[F^2]\right)^2\right]
		\right)^{\frac{1}{2}}+
		\left\lvert \mathbb{E}[F^2]-\sigma^2\right\rvert.
	\end{equation*}
	Finally, \eqref{eq:key-fixed-chaos} bounds the first term on the right-hand side by
	\begin{equation*}
		C_{p,q,q}
		\sqrt{\frac{q-1}{3q}}
		\sqrt{\kappa_4(F)}.
	\end{equation*}
	Substitution into the preceding Wasserstein estimate proves \eqref{eq:single-fourth}.
\end{proof}

Now let $q\neq r$, $F=I_q(f)$, $G=I_r(g)$, and $Z=F+G$.  For the sum of two Wiener-It\^o chaoses, Diez, Maia and Nourdin proved the following lemma with help of Chat-GPT-5.

\begin{lemma}\label{lem:two-chaos-variance}\cite{DMN25}
	With the preceding notation, if $q$ and $r$ have opposite parity, then the estimate
	\begin{equation}\label{eq:two-chaos-opposite}
		{\rm Var}(\Gamma_Z)
		\leq
		\frac{3}{2}\kappa_4(Z)
	\end{equation}
	holds.
\end{lemma}

Now we give the proof of Theorem~\ref{thm:two-chaos-Wp}.

\begin{proof}
	Let $\mu_\sigma$ be the distribution of $Z/\sigma$.  As in the proof of Theorem~\ref{thm:single-chaos-Wp}, scaling, the scalar Stein-kernel estimate, and conditional Jensen's inequality give
	\begin{equation}\label{eq:two-scaled-Stein}
		W_p\left(\mu,\gamma_\sigma\right)
		\leq
		\frac{C_p}{\sigma}
		\left(
		\mathbb{E}\left[\left\lvert \Gamma_Z-\sigma^2\right\rvert^p\right]
		\right)^{\frac{1}{p}}.
	\end{equation}
	More precisely,
	\begin{equation*}
		DZ=DF+DG,
		\qquad
		-DL^{-1}Z=\frac{1}{q}DF+\frac{1}{r}DG.
	\end{equation*}
	Thus both factors in $\Gamma_Z=\left\langle DZ,-DL^{-1}Z\right\rangle_{\mathbb{H}}$ contain only the chaos orders $q-1$ and $r-1$.  The product formula implies
	\begin{equation*}
		\Gamma_Z-\sigma^2
		\in
		\bigoplus_{\ell=0}^{2m-2}\mathcal H_\ell,
		\qquad m=q\vee r.
	\end{equation*}
	Applying \eqref{eq:finite-hypercontractivity} to \eqref{eq:two-scaled-Stein} yields
	\begin{equation*}
		W_p\left(\mu,\gamma_\sigma\right)
		\leq
		\frac{C_pB_{p,m,m}}{\sigma}
		\left(
		\mathbb{E}\left[\left(\Gamma_Z-\sigma^2\right)^2\right]
		\right)^{\frac{1}{2}},
	\end{equation*}
	which is \eqref{eq:two-first}.

	For the fourth-cumulant estimate, use $\mathbb{E}[\Gamma_Z]=\mathbb{E}[Z^2]$ and decompose
	\begin{equation*}
		\Gamma_Z-\sigma^2
		=
		\left(\Gamma_Z-\mathbb{E}[Z^2]\right)
		+
		\left(\mathbb{E}[Z^2]-\sigma^2\right).
	\end{equation*}
	The first term is centered and belongs to the sum of the first $2m-2$ nonconstant chaoses.  Hence Minkowski's inequality and \eqref{eq:centered-finite-hypercontractivity} give
	\begin{equation*}
		\left(
		\mathbb{E}\left[\left\lvert \Gamma_Z-\sigma^2\right\rvert^p\right]
		\right)^{\frac{1}{p}}
		\leq
		C_{p,m,m}
		\left(
		\mathbb{E}\left[\left(\Gamma_Z-\mathbb{E}[Z^2]\right)^2\right]
		\right)^{\frac{1}{2}}+
		\left\lvert \mathbb{E}[Z^2]-\sigma^2\right\rvert\\
	\end{equation*}
	Combining this inequality with \eqref{eq:two-scaled-Stein} and the opposite-parity estimate \eqref{eq:two-chaos-opposite} proves \eqref{eq:two-fourth}.
\end{proof}

\subsection{Proof of Theorem~\ref{thm:vector-chaos}}
Now we turn to the cases of vector-valued Wiener-It\^o chaos. For a random $d\times d$ matrix $A=(A_{ij})$, write
\begin{equation*}
	\left\lVert A\right\rVert_{p}
	:=
	\left(\sum_{i,j=1}^d\mathbb{E}\left\lvert A_{ij}\right\rvert^p\right)^{\frac{1}{p}}.
\end{equation*}
For a vector $F=(F_1,\ldots,F_d)$ of centered Wiener functionals, define
\begin{equation}\label{eq:Gamma-matrix}
	\Gamma_{ij}(F)
	:=
	\left\langle DF_j,-DL^{-1}F_i\right\rangle_{\mathbb{H}}.
\end{equation}
Then
\begin{equation*}
	\tau_F^{ij}(x)
	=
	\mathbb{E}\left[\Gamma_{ij}(F)\mid F=x\right]
\end{equation*}
defines a Stein kernel in the sense of \eqref{eq:vector-Stein-kernel}.

\begin{lemma}\label{lem:matrix-lp}
	For deterministic $d\times d$ matrices $A$ and random $d\times d$ matrices $M$,
	\begin{equation*}
		\left\lVert AMA^{\mathsf T}\right\rVert_p
		\leq
		d^{2\left\lvert \frac{1}{p}-\frac{1}{2}\right\rvert}\left\lVert A\right\rVert_{\mathrm{op}}^2\left\lVert M\right\rVert_p.
	\end{equation*}
\end{lemma}

\begin{proof}
	Set $D=d^2$ and identify a $d\times d$ matrix with a vector in $\mathbb{R}^{D}$ by stacking its columns.  The vectorization identity gives, pointwise on the underlying probability space, we have
	\begin{equation*}
		\operatorname{vec}\left(AMA^{\mathsf T}\right)
		=(A\otimes A)\operatorname{vec}(M),
	\end{equation*}
	where $\operatorname{vec}(M)$ is defined by $\operatorname{vec}: M_d(\mathbb{C})\longmapsto\mathbb{C}^{d\times d}$
	\begin{eqnarray*}	
		\operatorname{vec}(M)=(m_{11},m_{21},\cdots,m_{d1},m_{12},m_{22},\cdots,m_{d2},\cdots,m_{dd})^{T}.
	\end{eqnarray*}
	For every $x\in\mathbb{R}^{D}$, equivalence of the $\ell_p$- and $\ell_2$-norms yields
	\begin{equation*}
		\left\lVert x\right\rVert_{\ell_p}
		\leq
		D^{\left(\frac{1}{p}-\frac{1}{2}\right)_+}
		\left\lVert x\right\rVert_{\ell_2},
		\qquad
		\left\lVert x\right\rVert_{\ell_2}
		\leq
		D^{\left(\frac{1}{2}-\frac{1}{p}\right)_+}
		\left\lVert x\right\rVert_{\ell_p}.
	\end{equation*}
	Together with $\left\lVert A\otimes A\right\rVert_{\mathrm{op}}=\left\lVert A\right\rVert_{\mathrm{op}}^2$, these inequalities imply
	\begin{eqnarray*}
		\left\lVert \operatorname{vec}\left(AMA^{\mathsf T}\right)\right\rVert_{\ell_p}
		&\leq&
		D^{\left\lvert \frac{1}{p}-\frac{1}{2}\right\rvert}
		\left\lVert A\right\rVert_{\mathrm{op}}^2
		\left\lVert \operatorname{vec}(M)\right\rVert_{\ell_p}\\
		&=&
		d^{2\left\lvert \frac{1}{p}-\frac{1}{2}\right\rvert}
		\left\lVert A\right\rVert_{\mathrm{op}}^2
		\left\lVert \operatorname{vec}(M)\right\rVert_{\ell_p}.
	\end{eqnarray*}
	Raising both sides to the $p$th power and taking expectation gives the stated estimate.
\end{proof}

Notice that the exponent defined in \eqref{eq:beta} satisfies
\begin{equation*}
	\beta(p)
	=
	\alpha(p)+2\left\lvert \frac{1}{p}-\frac{1}{2}\right\rvert.
\end{equation*}

Now we give the proof of Theorem~\ref{thm:vector-chaos}.

\begin{proof}
	Let $Y=C^{-\frac{1}{2}}F$ and denote by $\mu_Y$ the distribution of $Y$.  Then
	\begin{equation*}
		W_p\left(\mu,\gamma_C\right)
		\leq
		\left\lVert C^{\frac{1}{2}}\right\rVert_{\mathrm{op}}
		W_p\left(\mu_Y,\gamma_d\right).
	\end{equation*}
	Indeed, if $(Y_0,G_0)$ is any coupling of $\mu_Y$ and $\gamma_d$, then $(C^{\frac{1}{2}}Y_0,C^{\frac{1}{2}}G_0)$ is a coupling of $\mu$ and $\gamma_C$, and
	\begin{equation*}
		\left|C^{\frac{1}{2}}Y_0-C^{\frac{1}{2}}G_0\right|
		\leq
		\left\lVert C^{\frac{1}{2}}\right\rVert_{\mathrm{op}}\left|Y_0-G_0\right|.
	\end{equation*}
	Taking the infimum over all couplings proves the preceding inequality.

	By the change of variables $F=C^{\frac{1}{2}}Y$, a Stein kernel of $Y$ is
	\begin{equation*}
		\tau_Y(y)
		=
		C^{-\frac{1}{2}}\tau_F\left(C^{\frac{1}{2}}y\right)C^{-\frac{1}{2}}.
	\end{equation*}
	Applying \eqref{eq:LNP-Wp} to $Y$, and then Lemma~\ref{lem:matrix-lp} with $A=C^{-\frac{1}{2}}$, gives
	\begin{eqnarray*}
		W_p\left(\mu,\gamma_C\right)
		&\leq&
		C_p d^{\alpha(p)}
		\left\lVert C^{\frac{1}{2}}\right\rVert_{\mathrm{op}}
		\left\lVert \tau_Y-I_d\right\rVert_p\\
		&\leq&
		C_p d^{\beta(p)}
		\left\lVert C^{\frac{1}{2}}\right\rVert_{\mathrm{op}}
		\left\lVert C^{-1}\right\rVert_{\mathrm{op}}
		\left\lVert \tau_F-C\right\rVert_p.
	\end{eqnarray*}
	Here the identity
	\begin{equation*}
		\beta(p)
		=
		\alpha(p)+2\left\lvert \frac{1}{p}-\frac{1}{2}\right\rvert
	\end{equation*}
	explains the dimension exponent.  Since
	\begin{equation*}
		\tau_F^{ij}(F)
		=
		\mathbb{E}\left[\Gamma_{ij}(F)\mid F\right],
	\end{equation*}
	conditional Jensen's inequality gives
	\begin{equation*}
		\left\lVert \tau_F-C\right\rVert_p^p
		\leq
		\sum_{i,j=1}^d
		\mathbb{E}\left[\left\lvert \Gamma_{ij}(F)-C_{ij}\right\rvert^p\right].
	\end{equation*}
	Finally, the product formula shows that $\Gamma_{ij}(F)-C_{ij}$ has chaos order at most $q_i+q_j-2$.  Lemma~\ref{lem:finite-hypercontractivity}, applied separately to every pair $(i,j)$, yields
	\begin{equation*}
		\mathbb{E}\left[\left\lvert \Gamma_{ij}(F)-C_{ij}\right\rvert^p\right]
		\leq
		B_{p,q_i,q_j}^p
		\mathbb{E}\left[\left\lvert \Gamma_{ij}(F)-C_{ij}\right\rvert^2\right]^{\frac{p}{2}}.
	\end{equation*}
	Substituting this estimate into the preceding display proves \eqref{eq:vector-exact-first}.

	For the second assertion, the standard contraction estimates of \cite{NPR10} imply that, for $i\neq j$ and a constant $E_{q,r}$ depending only on the two chaos orders,
\begin{equation}
	\mathbb{E}\left[
	\left\lvert \Gamma_{ij}(F)-\mathbb{E}[F_iF_j]\right\rvert^2
	\right]
	\leq
	E_{q_i,q_j}
	\bigg(
	\mathbb{E}[F_i^2]\sqrt{\kappa_4(F_j)}
	+
	\mathbb{E}[F_j^2]\sqrt{\kappa_4(F_i)}
	+
	\kappa_4(F_i)+\kappa_4(F_j)
	\bigg),
	\label{eq:cross-contraction-bound}
\end{equation}
where 
	\begin{equation*}
	E_{q,r}=\max\left\{2,\left(\sum_{s=1}^{q-1}\binom{2s}{s}\right),\frac{1}{2}\sum_{s=1}^{(q\wedge r)-1}(q+r-2s)!\binom{r}{s}^2,\frac{1}{2}\sum_{s=1}^{(q\wedge r)-1}(q+r-2s)!\binom{q}{s}^2 \right\}.
\end{equation*}
	By the triangle inequality,
\begin{eqnarray*}
	\left\lVert \Gamma_{ij}(F)-C_{ij}\right\rVert_{L^2}
	&\leq&
	\left\lVert \Gamma_{ij}(F)-\mathbb{E}[F_iF_j]\right\rVert_{L^2}
	+
	\left\lvert \mathbb{E}[F_iF_j]-C_{ij}\right\rvert.
\end{eqnarray*}
	For $i=j$, estimate \eqref{eq:key-fixed-chaos} gives
\begin{equation*}
	\left\lVert \Gamma_{ii}(F)-\mathbb{E}[F_i^2]\right\rVert_{L^2}
	\leq
	\sqrt{\frac{q_i-1}{3q_i}}\sqrt{\kappa_4(F_i)}.
\end{equation*}
	For $i\neq j$, take square roots in \eqref{eq:cross-contraction-bound} and use the elementary inequality $\sqrt{x_1+\cdots+x_J}\leq\sqrt{x_1}+\cdots+\sqrt{x_J}$ for $x_1,\ldots,x_J\geq0$. We obtain
\begin{equation*}
	\left\lVert \Gamma_{ij}(F)-\mathbb{E}[F_iF_j]\right\rVert_{L^2}
	\leq
	C_{q_i,q_j}
	\bigg(
	\mathbb{E}[F_i^2]^{\frac{1}{2}}\kappa_4(F_j)^{\frac{1}{4}}
	+
	\mathbb{E}[F_j^2]^{\frac{1}{2}}\kappa_4(F_i)^{\frac{1}{4}}
	+
	\sqrt{\kappa_4(F_i)}
	+
	\sqrt{\kappa_4(F_j)}
	\bigg).
\end{equation*}
	Substitute these estimates into \eqref{eq:vector-exact-first} and sum over $i,j$.  Absorbing the fixed factors $C_p$, $d^{\beta(p)}$, $B_{p,q_i,q_j}$, and $E_{q_i,q_j}$ into $D_{p,d,q_1,\ldots,q_d}$ proves \eqref{eq:vector-fourth-intro}.
\end{proof}

For completeness, the same argument covers vectors whose components are sums of two chaoses of opposite parity. Let
\begin{equation*}
	Z_i=F_i+G_i,
	\qquad
	F_i=I_{q_i}(f_i),
	\qquad
	G_i=I_{r_i}(g_i),
	\qquad
	q_i\neq r_i,
\end{equation*}
where $q_i$ and $r_i$ have opposite parity, and put $m_i=q_i\vee r_i$. Also set
\begin{equation*}
	v_i:=\mathbb{E}[Z_i^2]=\mathbb{E}[F_i^2]+\mathbb{E}[G_i^2].
\end{equation*}

\subsection{Proof of Theorem~\ref{thm:vector-two-chaos}}

\begin{proof}
	Apply the linear transformation $C^{-\frac{1}{2}}$ exactly as in the proof of Theorem~\ref{thm:vector-chaos}.  The vector Stein-kernel estimate, Lemma~\ref{lem:matrix-lp}, and conditional Jensen's inequality give
	\begin{equation*}
		W_p\left(\mu,\gamma_C\right)
		\leq
		C_p d^{\beta(p)}
		\left\lVert C^{\frac{1}{2}}\right\rVert_{\mathrm{op}}
		\left\lVert C^{-1}\right\rVert_{\mathrm{op}}
		\left(
		\sum_{i,j=1}^d
		\mathbb{E}\left[\left\lvert \Gamma_{ij}(Z)-C_{ij}\right\rvert^p\right]
		\right)^{\frac{1}{p}}.
	\end{equation*}
	Since $Z_i$ is the sum of chaoses of orders $q_i$ and $r_i$, the product formula shows that $\Gamma_{ij}(Z)-C_{ij}$ belongs to
	\begin{equation*}
		\bigoplus_{\ell=0}^{m_i+m_j-2}\mathcal H_\ell.
	\end{equation*}
	Therefore, Lemma~\ref{lem:finite-hypercontractivity} yields
	\begin{equation*}
		\mathbb{E}\left[\left\lvert \Gamma_{ij}(Z)-C_{ij}\right\rvert^p\right]
		\leq
		B_{p,m_i,m_j}^p
		\mathbb{E}\left[\left\lvert \Gamma_{ij}(Z)-C_{ij}\right\rvert^2\right]^{\frac{p}{2}}.
	\end{equation*}
	Substituting this estimate entrywise proves \eqref{eq:vector-two-exact}.

	For \eqref{eq:vector-two-fourth}, define, for centered Wiener functionals $U$ and $V$,
	\begin{equation*}
		\Gamma(U,V)
		:=
		\left\langle DV,-DL^{-1}U\right\rangle_{\mathbb{H}}.
	\end{equation*}
	Then $\mathbb{E}[\Gamma(U,V)]=\mathbb{E}[UV]$, and
	\begin{eqnarray*}
		\Gamma_{ij}(Z)
		&=&
		\Gamma(F_i,F_j)+\Gamma(F_i,G_j)+
		\Gamma(G_i,F_j)+\Gamma(G_i,G_j).
	\end{eqnarray*}
	Thus the centered quantity $\Gamma_{ij}(Z)-\mathbb{E}[Z_iZ_j]$ is the sum of four centered pairings.  For $i\neq j$, apply the contraction estimate \eqref{eq:cross-contraction-bound} to each pair
	\begin{equation*}
		(U,V)
		\in
		\left\{F_i,G_i\right\}\times\left\{F_j,G_j\right\}.
	\end{equation*}
	After taking square roots, each such term is bounded by a constant depending only on the corresponding chaos orders times
	\begin{equation*}
		\mathbb{E}[U^2]^{\frac{1}{2}}\kappa_4(V)^{\frac{1}{4}}
		+
		\mathbb{E}[V^2]^{\frac{1}{2}}\kappa_4(U)^{\frac{1}{4}}
		+
		\kappa_4(U)^{\frac{1}{2}}
		+
		\kappa_4(V)^{\frac{1}{2}}.
	\end{equation*}
	The fourth-cumulant decomposition used in Lemma~\ref{lem:two-chaos-variance} gives
	\begin{equation*}
		\kappa_4(F_i)+\kappa_4(G_i)
		\leq
		\kappa_4(Z_i),
	\end{equation*}
	while $\mathbb{E}[F_i^2]\leq v_i$ and $\mathbb{E}[G_i^2]\leq v_i$.  It follows in particular that
	\begin{equation*}
		\kappa_4(F_i)^{\frac{1}{4}}
		+
		\kappa_4(G_i)^{\frac{1}{4}}
		\leq
		2^{\frac{3}{4}}\kappa_4(Z_i)^{\frac{1}{4}},
	\end{equation*}
	and the analogous estimate holds with the exponent $\dfrac{1}{2}$.  The triangle inequality in $L^2$, followed by summation over the four pairings, therefore yields
	\begin{equation*}
		\left\lVert \Gamma_{ij}(Z)-\mathbb{E}[Z_iZ_j]\right\rVert_{L^2}
		\leq
		C_{\boldsymbol q,\boldsymbol r}
		\left(
		v_i^{\frac{1}{2}}\kappa_4(Z_j)^{\frac{1}{4}}
		+
		v_j^{\frac{1}{2}}\kappa_4(Z_i)^{\frac{1}{4}}
		+
		\kappa_4(Z_i)^{\frac{1}{2}}
		+
		\kappa_4(Z_j)^{\frac{1}{2}}
		\right).
	\end{equation*}
	For $i=j$, Lemma~\ref{lem:two-chaos-variance} gives the direct estimate
	\begin{equation*}
		\left(
		\mathbb{E}\left[\left(\Gamma_{ii}(Z)-\mathbb{E}[Z_i^2]\right)^2\right]
		\right)^{\frac{1}{2}}
		\leq
		\sqrt{2\kappa_4(Z_i)}.
	\end{equation*}
	Finally, for every $i,j$,
	\begin{equation*}
		\left\lVert \Gamma_{ij}(Z)-C_{ij}\right\rVert_{L^2}
		\leq
		\left\lVert \Gamma_{ij}(Z)-\mathbb{E}[Z_iZ_j]\right\rVert_{L^2}
		+
		\left\lvert \mathbb{E}[Z_iZ_j]-C_{ij}\right\rvert.
	\end{equation*}
	Substitution into \eqref{eq:vector-two-exact}, followed by summation over $i,j$, proves \eqref{eq:vector-two-fourth} after all order-dependent constants are absorbed into $D_{p,d,\boldsymbol q,\boldsymbol r}$.
\end{proof}

\begin{remark}
	When $p=2$, one has $C_2=1$ and $\beta(2)=0$.  Thus \eqref{eq:vector-exact-first} is dimension-free apart from the entrywise matrix norm.  Since $W_1\leq W_2$, it also recovers the usual multivariate $W_1$ estimate of \cite{NPR10}, up to the chosen convention for the matrix norm.
\end{remark}
\subsection{Proof of Theorem~\ref{thm:finite-chaos-sum}}

The $d_2$ estimate holds for every centered $F\in\mathbb D^{1,2}$ with unit variance. The finite-chaos assumption is used only to obtain the explicit contraction expansion and the $L^p$ estimate.

\begin{proof}
	Let $N$ be a standard Gaussian random variable with distribution $\gamma$, independent of the Gaussian process generating $F$.  The normalization $h(0)=h'(0)=0$ causes no loss, because $F$ and $N$ have the same mean.  Put
	\begin{equation*}
		F_t=\sqrt{t}\,F+\sqrt{1-t}\,N,
		\qquad 0\leq t\leq1.
	\end{equation*}
	For $h\in C^2(\mathbb{R})$ with bounded second derivative, define $\Psi(t)=\mathbb{E}[h(F_t)]$. Integration by parts formula gives
	\begin{eqnarray*}
		\mathbb{E}\left[Fh'(F_t)\right]
		&=&
		\sqrt{t}\,\mathbb{E}\left[h''(F_t)\Gamma_F\right],\\
		\mathbb{E}\left[Nh'(F_t)\right]
		&=&
		\sqrt{1-t}\,\mathbb{E}\left[h''(F_t)\right].
	\end{eqnarray*}
	Consequently,
	\begin{equation*}
		\Psi'(t)
		=
		\frac{1}{2}\mathbb{E}\left[h''(F_t)\left(\Gamma_F-1\right)\right].
	\end{equation*}
	Integration over $t\in[0,1]$ proves \eqref{eq:smart-path-defect}. 
	Write $F_q=I_q(f_q)$.  By \eqref{eq:D-chaos} and the product formula,
	\begin{equation*}
		\left\langle DF_q,DF_r\right\rangle_{\mathbb{H}}
		=
		qr\sum_{k=1}^{q\wedge r}
		(k-1)!\binom{q-1}{k-1}\binom{r-1}{k-1}
		I_{q+r-2k}\left(f_q\widetilde{\otimes}_k f_r\right).
	\end{equation*}
	Since
	\begin{equation*}
		-DL^{-1}F=\sum_{r=1}^{M}\frac{1}{r}DF_r,
	\end{equation*}
	we obtain
	\begin{equation*}
		\Gamma_F
		=
		\sum_{q,r=1}^{M}\sum_{k=1}^{q\wedge r}
		q(k-1)!\binom{q-1}{k-1}\binom{r-1}{k-1}
		I_{q+r-2k}\left(f_q\widetilde{\otimes}_k f_r\right).
	\end{equation*}
	The zeroth-chaos term occurs only when $q=r=k$ and equals
	\begin{equation*}
		\sum_{q=1}^{M}q!\left\lVert f_q\right\rVert_{\mathbb{H}^{\otimes q}}^2
		=
		\mathbb{E}[F^2]
		=1.
	\end{equation*}
	Thus
	\begin{equation*}
		\Gamma_F-1=\sum_{\ell=1}^{2M-2}I_\ell(h_\ell).
	\end{equation*}
	Orthogonality proves \eqref{eq:finite-sum-exact-variance}. Finally, \eqref{eq:LNP-Wp}, conditional Jensen, and \eqref{eq:centered-finite-hypercontractivity} give
	\begin{equation*}
		W_p\left(\mu,\gamma\right)
		\leq
		C_p\left\lVert \Gamma_F-1\right\rVert_{L^p}
		\leq
		C_pC_{p,M,M}
		\left(
		\mathbb{E}\left[\left(\Gamma_F-1\right)^2\right]
		\right)^{\frac{1}{2}},
	\end{equation*}
	which is \eqref{eq:finite-sum-Wp}.
\end{proof}

\begin{remark}\label{rem:no-general-fourth-cumulant}
	Formula \eqref{eq:finite-sum-exact-variance} is the desired explicit estimate for a finite sum of Wiener--It\^o chaoses.  It also explains why a bound involving only $\kappa_4(F)$ is unavailable for an arbitrary chaos sum: mixed contractions $f_q\widetilde{\otimes}_k f_r$ enter the Malliavin defect, and different pairs may contribute to the same chaos level.  Additional structure is therefore needed.  Such structure is present for a fixed chaos, for independent chaos components, and for a sum of two chaoses of opposite parity.  Without it, a fourth moment theorem can fail; see the same-parity examples in \cite{BKS25}.
\end{remark}

\section{Diagonal Wiener--It\^o chaos of arbitrary order}\label{sec:diagonal-chaos}

In this section, $n\geq2$ is fixed.  We derive a family of estimates without using Stein's method.

\subsection{Diagonal representation and Hermite derivatives}

Let $f\in\mathbb{H}^{\odot n}$ be
\begin{equation}\label{eq:diagonal-kernel}
	f=\sum_{k\geq1}\lambda_k e_k^{\otimes n}
\end{equation}
in $\mathbb{H}^{\otimes n}$, where $(e_k)_{k\geq1}$ is orthonormal and $(\lambda_k)_{k\geq1}\in\ell^2$.  Such a tensor will be called \emph{orthogonally diagonalizable}.  If $n=2$, every symmetric kernel has this representation after identifying it with a compact self-adjoint Hilbert--Schmidt operator.  If $n\geq3$, \eqref{eq:diagonal-kernel} is an additional hypothesis and is not a spectral theorem for arbitrary symmetric kernels.

Set
\begin{equation*}
	\xi_k:=X(e_k),
	\qquad k\geq1.
\end{equation*}
The random variables $(\xi_k)_{k\geq1}$ are independent and standard Gaussian.  Since
\begin{equation*}
	I_n\left(e_k^{\otimes n}\right)=H_n(\xi_k),
\end{equation*}
linearity and the Wiener--It\^o isometry yield
\begin{equation}\label{eq:diagonal-chaos}
	F:=I_n(f)
	=
	\sum_{k\geq1}\lambda_k H_n(\xi_k)
\end{equation}
in $L^2(\Omega)$.  Indeed, for $L>M$,
\begin{equation*}
	\mathbb{E}\left[
	\left(
	\sum_{k=M+1}^{L}\lambda_kH_n(\xi_k)
	\right)^2
	\right]
	=
	n!\sum_{k=M+1}^{L}\lambda_k^2
	\longrightarrow0.
\end{equation*}
We impose the normalization
\begin{equation}\label{eq:diagonal-normalization}
	\mathbb{E}[F^2]=n!\sum_{k\geq1}\lambda_k^2=1.
\end{equation}

Define Hermite polynomials by
\begin{equation*}
	H_m(x)=(-1)^m e^{\frac{x^2}{2}}\frac{d^m}{dx^m}e^{-\frac{x^2}{2}}.
\end{equation*}
One has
\begin{equation}\label{eq:Hermite-derivative}
	H_m'(x)=mH_{m-1}(x),
\end{equation}
as well as the three-term recurrence
\begin{equation}\label{eq:Hermite-three-term}
	H_{m+1}(x)
	=
	xH_m(x)-mH_{m-1}(x),
	\qquad
	H_0(x)=1,
	\quad
	H_1(x)=x.
\end{equation}
More generally,
\begin{equation}\label{eq:Hermite-higher-derivative}
	\frac{d^j}{dx^j}H_m(x)
	=
	\frac{m!}{(m-j)!}H_{m-j}(x),
	\qquad 0\leq j\leq m.
\end{equation}
Since $D\xi_k=e_k$, the chain rule gives
\begin{equation}\label{eq:DjF-diagonal}
	D^jF
	=
	\frac{n!}{(n-j)!}
	\sum_{k\geq1}\lambda_kH_{n-j}(\xi_k)e_k^{\otimes j},
	\qquad 1\leq j\leq n.
\end{equation}
In particular,
\begin{eqnarray}
	DF
	&=&
	n\sum_{k\geq1}\lambda_kH_{n-1}(\xi_k)e_k,
	\label{eq:DF-diagonal}\\
	D^2F
	&=&
	n(n-1)\sum_{k\geq1}\lambda_kH_{n-2}(\xi_k)e_k^{\otimes2}.
	\label{eq:D2F-diagonal}
\end{eqnarray}
Thus the Malliavin covariance of $F$ is
\begin{equation}\label{eq:gamma-diagonal}
	\gamma_F:=\left\lVert DF\right\rVert_{\mathbb{H}}^2
	=
	n^2\sum_{k\geq1}\lambda_k^2H_{n-1}(\xi_k)^2.
\end{equation}

\subsection{Hermite small-ball estimates and negative moments}
In this subsection, we follows the idea in \cite{IS89,IS04} to derive the conditions to guarantee that Wiener chaos $F$ with the form \ref {eq:diagonal-chaos} are smooth in the sense of Malliavin calculus.
\begin{lemma}\label{lem:Hermite-small-ball}
	Let $m\geq1$ and $\xi\sim\mathcal{N}(0,1)$.  There exists a constant $C_m>0$ such that
	\begin{equation}\label{eq:Hermite-small-ball}
		\mathbb{P}\left(\left\lvert H_m(\xi)\right\rvert\leq u\right)
		\leq
		C_m u,
		\qquad 0<u\leq1.
	\end{equation}
\end{lemma}

\begin{proof}
	Let
	\begin{equation*}
		a_1<a_2<\cdots<a_m
	\end{equation*}
	be the zeros of $H_m$.  Since the Hermite polynomials are orthogonal with respect to the strictly positive Gaussian weight, the standard zero theorem for orthogonal polynomials shows that all $m$ zeros are real.  They are also simple.  Indeed, if $H_m(a)=H_m'(a)=0$, then \eqref{eq:Hermite-derivative} gives $H_{m-1}(a)=0$.  Repeated use of \eqref{eq:Hermite-three-term} then gives
	\begin{equation*}
		H_{m-2}(a)=\cdots=H_0(a)=0,
	\end{equation*}
	contradicting $H_0(a)=1$.  Consequently, by \eqref{eq:Hermite-derivative},
	\begin{equation*}
		H_m'(a_j)=mH_{m-1}(a_j)\neq0.
	\end{equation*}
	For each $j$, continuity of $H_m'$ yields numbers $\delta_j,c_j>0$ such that
	\begin{equation*}
		\left\lvert x-a_j\right\rvert\leq\delta_j
		\quad\Longrightarrow\quad
		\left\lvert H_m(x)\right\rvert\geq c_j\left\lvert x-a_j\right\rvert.
	\end{equation*}
	After decreasing the $\delta_j$ if necessary, the intervals
	\begin{equation*}
		(a_j-\delta_j,a_j+\delta_j)
	\end{equation*}
	are disjoint.  On the closed set outside their union, $H_m$ has no zero and tends to infinity in absolute value as $\left\lvert x\right\rvert\to\infty$.  Hence there is $\eta_m>0$ such that
	\begin{equation*}
		\left\lvert H_m(x)\right\rvert\geq\eta_m
	\end{equation*}
	on that set.  Put $c=\min_jc_j$.  If $0<u<\eta_m$, then
	\begin{equation*}
		\left\{x:\left\lvert H_m(x)\right\rvert\leq u\right\}
		\subset
		\bigcup_{j=1}^m
		\left[a_j-\frac{u}{c},a_j+\frac{u}{c}\right].
	\end{equation*}
	Therefore, with $\phi$ denoting the standard Gaussian density,
	\begin{eqnarray*}
		\mathbb{P}\left(\left\lvert H_m(\xi)\right\rvert\leq u\right)
		&\leq&
		\sum_{j=1}^m
		\int_{a_j-\frac{u}{c}}^{a_j+\frac{u}{c}}\phi(x)\,dx\\
		&\leq&
		\frac{2m}{c\sqrt{2\pi}}u.
	\end{eqnarray*}
	For $\eta_m\leq u\leq1$, the estimate follows from
	\begin{equation*}
		\mathbb{P}\left(\left\lvert H_m(\xi)\right\rvert\leq u\right)
		\leq1\leq\eta_m^{-1}u.
	\end{equation*}
\end{proof}

\begin{proposition}\label{prop:negative-moments}
	Choose $r$ nonzero coefficients in \eqref{eq:diagonal-kernel}, relabel them as $\lambda_1,\ldots,\lambda_r$, and put
	\begin{equation*}
		\lambda_{*,r}:=\min_{1\leq k\leq r}\left\lvert \lambda_k\right\rvert>0.
	\end{equation*}
	There exists $C_{n,r}>0$ such that, for all sufficiently small $\varepsilon>0$,
	\begin{equation}\label{eq:gamma-small-ball}
		\mathbb{P}\left(\gamma_F\leq\varepsilon\right)
		\leq
		C_{n,r}\lambda_{*,r}^{-r}\varepsilon^{\frac{r}{2}}.
	\end{equation}
	Consequently,
	\begin{equation}\label{eq:finite-negative-moment}
		\mathbb{E}\left[\gamma_F^{-p}\right]<\infty
		\qquad\text{for every }0<p<\frac{r}{2}.
	\end{equation}
\end{proposition}

\begin{proof}
	By \eqref{eq:gamma-diagonal},
	\begin{equation*}
		\gamma_F
		\geq
		n^2\lambda_{*,r}^2
		\sum_{k=1}^rH_{n-1}(\xi_k)^2.
	\end{equation*}
	Hence
	\begin{equation*}
		\left\{\gamma_F\leq\varepsilon\right\}
		\subset
		\bigcap_{k=1}^r
		\left\{
		\left\lvert H_{n-1}(\xi_k)\right\rvert
		\leq
		\frac{\sqrt\varepsilon}{n\lambda_{*,r}}
		\right\}.
	\end{equation*}
	The variables $\xi_1,\ldots,\xi_r$ are independent.  Lemma~\ref{lem:Hermite-small-ball} therefore gives
	\begin{equation*}
		\mathbb{P}\left(\gamma_F\leq\varepsilon\right)
		\leq
		\left(
		C_n\frac{\sqrt\varepsilon}{n\lambda_{*,r}}
		\right)^r,
	\end{equation*}
	which proves \eqref{eq:gamma-small-ball}.  For a positive random variable $Y$ and $p>0$,
	\begin{equation}\label{eq:negative-layer-cake}
		\mathbb{E}\left[Y^{-p}\right]
		=
		p\int_0^\infty
		u^{-p-1}\mathbb{P}(Y<u)\,du.
	\end{equation}
	The integral over $[1,\infty)$ is finite. Near zero, estimate \eqref{eq:gamma-small-ball} makes the integrand bounded by a constant times
	\begin{equation*}
		u^{\frac{r}{2}-p-1},
	\end{equation*}
	which is integrable at zero exactly when $p<\dfrac{r}{2}$; on every compact interval bounded away from zero, the integral is finite.
\end{proof}

\begin{corollary}\label{cor:smooth-nondegenerate}
	If infinitely many coefficients $\lambda_k$ are nonzero, then
	\begin{equation*}
		F\in\mathbb D^\infty,
		\qquad
		\gamma_F^{-1}\in\bigcap_{p>0}L^p(\Omega).
	\end{equation*}
	Thus $F$ is smoothly nondegenerate in the sense of Malliavin calculus.
\end{corollary}

\begin{proof}
	Fix $0\leq j\leq n$.  Formula \eqref{eq:DjF-diagonal} and the orthogonality of $(e_k)_{k\geq1}$ give
	\begin{equation*}
		\mathbb{E}\left[\left\lVert D^jF\right\rVert_{\mathbb{H}^{\otimes j}}^2\right]
		=
		\left(\frac{n!}{(n-j)!}\right)^2
		(n-j)!\sum_{k\geq1}\lambda_k^2
		<\infty.
	\end{equation*}
	Thus the series defining $D^jF$ converges in $L^2\left(\Omega;\mathbb{H}^{\otimes j}\right)$.  Since it belongs to the $(n-j)$th Hilbert-valued Wiener--It\^o chaos, Hilbert-valued hypercontractivity implies
	\begin{equation*}
		\mathbb{E}\left[\left\lVert D^jF\right\rVert_{\mathbb{H}^{\otimes j}}^q\right]
		<\infty
		\qquad\text{for every }q>1.
	\end{equation*}
	For $j>n$, one has $D^jF=0$.  It follows that $F\in\mathbb D^{j,q}$ for every $j\geq0$ and every $q>1$, and hence $F\in\mathbb D^\infty$.

	It remains to verify the nondegeneracy of the Malliavin covariance.  Fix an arbitrary $p>0$.  Since infinitely many coefficients are nonzero, one may choose an integer $r>2p$ and then select $r$ nonzero coefficients.  For this finite selection, the number $\lambda_{*,r}$ in Proposition~\ref{prop:negative-moments} is strictly positive.  Since $p<\dfrac{r}{2}$, that proposition yields
	\begin{equation*}
		\mathbb{E}\left[\gamma_F^{-p}\right]<\infty.
	\end{equation*}
	In particular, $\mathbb{P}(\gamma_F=0)=0$.  As $p>0$ was arbitrary, $\gamma_F^{-1}$ belongs to $L^p(\Omega)$ for every $p>0$.  Combining this fact with $F\in\mathbb D^\infty$ proves smooth nondegeneracy.
\end{proof}

\begin{remark}
	For $n=2$, if the kernel has finite rank $R$, then $\gamma_F$ is a weighted chi-square random variable and
	\begin{equation*}
		\mathbb{E}\left[\gamma_F^{-p}\right]<\infty
		\quad\Longleftrightarrow\quad
		p<\frac{R}{2}.
	\end{equation*}
	For $n\geq3$, Proposition~\ref{prop:negative-moments} gives the same sufficient threshold from the number of nonzero diagonal coefficients.  It is not asserted to be necessary, because the local behavior of $H_{n-1}$ at its several zeros is more complicated than the chi-square case.
\end{remark}

\subsection{Fourth cumulants and a Berry--Esseen estimate}

Let $\xi\sim\mathcal{N}(0,1)$ and define
\begin{equation}\label{eq:cn}
	c_n:=\kappa_4\left(H_n(\xi)\right)
	=
	\mathbb{E}\left[H_n(\xi)^4\right]-3(n!)^2.
\end{equation}
The Hermite product formula
\begin{equation}\label{eq:Hermite-product}
	H_n(x)^2
	=
	\sum_{k=0}^n
	k!\binom{n}{k}^2H_{2n-2k}(x)
\end{equation}
and orthogonality imply
\begin{equation}\label{eq:Hermite-fourth}
	\mathbb{E}\left[H_n(\xi)^4\right]
	=
	\sum_{k=0}^n
	\left(k!\binom{n}{k}^2\right)^2(2n-2k)!.
\end{equation}
In particular, $c_n>0$ for $n\geq2$.  Indeed, retaining only the terms $k=0$ and $k=n$ gives
\begin{equation*}
	c_n
	\geq
	(2n)!-2(n!)^2
	=
	\left(\binom{2n}{n}-2\right)(n!)^2>0.
\end{equation*}

\begin{proposition}\label{prop:diagonal-fourth}
	For $F$ in \eqref{eq:diagonal-chaos},
	\begin{equation}\label{eq:diagonal-kappa}
		\kappa_4(F)
		=
		c_n\sum_{k\geq1}\lambda_k^4.
	\end{equation}
\end{proposition}

\begin{proof}
	For the finite truncation
	\begin{equation*}
		F^{(M)}=\sum_{k=1}^{M}\lambda_kH_n(\xi_k),
	\end{equation*}
	the summands $Y_k:=\lambda_kH_n(\xi_k)$ are independent and centered.  Expanding the fourth moment and using independence gives
	\begin{eqnarray*}
		\mathbb{E}\left[\left(F^{(M)}\right)^4\right]
		&=&
		\sum_{k=1}^{M}\mathbb{E}\left[Y_k^4\right]
		+
		6\sum_{1\leq k<\ell\leq M}
		\mathbb{E}\left[Y_k^2\right]\mathbb{E}\left[Y_\ell^2\right],\\
		\mathbb{E}\left[\left(F^{(M)}\right)^2\right]^2
		&=&
		\sum_{k=1}^{M}\mathbb{E}\left[Y_k^2\right]^2
		+
		2\sum_{1\leq k<\ell\leq M}
		\mathbb{E}\left[Y_k^2\right]\mathbb{E}\left[Y_\ell^2\right].
	\end{eqnarray*}
	Therefore the mixed terms cancel in the fourth cumulant, and
	\begin{eqnarray*}
		\kappa_4\left(F^{(M)}\right)
		&=&
		\sum_{k=1}^{M}
		\left(
		\mathbb{E}\left[Y_k^4\right]
		-3\mathbb{E}\left[Y_k^2\right]^2
		\right)\\
		&=&
		c_n\sum_{k=1}^{M}\lambda_k^4.
	\end{eqnarray*}
	This is the additivity of fourth cumulants written explicitly in the present setting.

	To pass to the infinite sum, the isometry gives
	\begin{equation*}
		\mathbb{E}\left[\left(F-F^{(M)}\right)^2\right]
		=
		n!\sum_{k>M}\lambda_k^2
		\longrightarrow0.
	\end{equation*}
	Moreover, $F-F^{(M)}$ belongs to the $n$th Wiener--It\^o chaos.  Hypercontractivity with exponent $4$ therefore yields
	\begin{equation*}
		\left(
		\mathbb{E}\left[\left\lvert F-F^{(M)}\right\rvert^4\right]
		\right)^{\frac{1}{4}}
		\leq
		3^{\frac{n}{2}}
		\left(
		\mathbb{E}\left[\left(F-F^{(M)}\right)^2\right]
		\right)^{\frac{1}{2}}
		\longrightarrow0.
	\end{equation*}
	Hence both the second and fourth moments of $F^{(M)}$ converge to the corresponding moments of $F$, so $\kappa_4\left(F^{(M)}\right)\to\kappa_4(F)$.  Finally,
	\begin{equation*}
		\sum_{k\geq1}\lambda_k^4
		\leq
		\left(\sum_{k\geq1}\lambda_k^2\right)^2
		<\infty,
	\end{equation*}
	and therefore $\sum_{k=1}^{M}\lambda_k^4\to\sum_{k\geq1}\lambda_k^4$.  Letting $M\to\infty$ proves \eqref{eq:diagonal-kappa}.
\end{proof}

\begin{proposition}\label{prop:diagonal-BE}
	Denote by $\mu_F$ the distribution of $F$ and by $\gamma$ the standard Gaussian distribution.  There exists a constant $C_n>0$, depending only on $n$, such that
	\begin{equation}\label{eq:diagonal-Kolmogorov}
		d_{\mathrm K}\left(\mu_F,\gamma\right)
		\leq
		C_n\sqrt{\kappa_4(F)}.
	\end{equation}
\end{proposition}

\begin{proof}
	Put
	\begin{equation*}
		Y_k=\lambda_kH_n(\xi_k),
		\qquad
		F^{(M)}=\sum_{k=1}^{M}Y_k,
		\qquad
		\sigma_M^2=n!\sum_{k=1}^{M}\lambda_k^2.
	\end{equation*}
	For fixed $F$, the integer $M$ is only a truncation parameter.  By \eqref{eq:diagonal-normalization}, $\sigma_M^2\to1$ and
	\begin{equation*}
		S_M:=\frac{F^{(M)}}{\sigma_M}\longrightarrow F
		\qquad\text{in }L^2.
	\end{equation*}
	Denote by $\mu_M$ the distribution of $S_M$.
	Here $S_M$ is defined for all sufficiently large $M$, when $\sigma_M>0$.  More explicitly,
	\begin{eqnarray*}
		\left\lVert S_M-F\right\rVert_{L^2}
		&\leq&
		\left\lvert \sigma_M^{-1}-1\right\rvert\left\lVert F^{(M)}\right\rVert_{L^2}
		+
		\left\lVert F^{(M)}-F\right\rVert_{L^2}\\
		&=&
		\left\lvert 1-\sigma_M\right\rvert
		+
		\left\lVert F^{(M)}-F\right\rVert_{L^2}
		\longrightarrow0.
	\end{eqnarray*}
	The classical Berry--Esseen theorem \cite{PET75} applied to the finite independent sum gives
	\begin{eqnarray}
		d_{\mathrm K}\left(\mu_M,\gamma\right)
		&\leq&
		C\sigma_M^{-3}
		\sum_{k=1}^{M}\mathbb{E}\left\lvert Y_k\right\rvert^3
		\label{eq:BE-finite}\\
		&=&
		C\sigma_M^{-3}\mathbb{E}\left\lvert H_n(\xi)\right\rvert^3
		\sum_{k=1}^{M}\left\lvert \lambda_k\right\rvert^3.
	\nonumber
	\end{eqnarray}
	Here and below, $C$ denotes a universal numerical constant whose value may change from line to line; in particular, it is independent of $M$ and of the coefficients $\lambda_k$.
	To remove the truncation, fix $\varepsilon>0$.  For every $x\in\mathbb{R}$,
	\begin{eqnarray*}
		\mathbb{P}(F\leq x)
		&\leq&
		\mathbb{P}(S_M\leq x+\varepsilon)
		+
		\mathbb{P}\left(\left\lvert F-S_M\right\rvert>\varepsilon\right),\\
		\mathbb{P}(S_M\leq x-\varepsilon)
		&\leq&
		\mathbb{P}(F\leq x)
		+
		\mathbb{P}\left(\left\lvert F-S_M\right\rvert>\varepsilon\right).
	\end{eqnarray*}
	It follows that
	\begin{equation*}
		d_{\mathrm K}\left(\mu_F,\gamma\right)
		\leq
		d_{\mathrm K}\left(\mu_M,\gamma\right)
		+
		\sup_{x\in\mathbb{R}}\left\lvert \Phi(x+\varepsilon)-\Phi(x)\right\rvert+
		\mathbb{P}\left(\left\lvert F-S_M\right\rvert>\varepsilon\right).
	\end{equation*}
	Let $M\to\infty$ and then let $\varepsilon\downarrow0$, note that $\Phi$ is uniformly continuous, we have
	\begin{equation*}
		d_{\mathrm K}\left(\mu_F,\gamma\right)
		\leq
		C_n\sum_{k\geq1}\left\lvert \lambda_k\right\rvert^3.
	\end{equation*}
	The series is finite because $(\lambda_k)\in\ell^2$.  Finally, Cauchy--Schwarz, \eqref{eq:diagonal-normalization}, and \eqref{eq:diagonal-kappa} yield
	\begin{equation*}
		\sum_{k\geq1}\left\lvert \lambda_k\right\rvert^3
		\leq
		\left(\sum_{k\geq1}\lambda_k^2\right)^{\frac{1}{2}}
		\left(\sum_{k\geq1}\lambda_k^4\right)^{\frac{1}{2}}
		=
		\frac{1}{\sqrt{n!c_n}}\sqrt{\kappa_4(F)}.
	\end{equation*}
\end{proof}

\subsection{Uniform inverse moments along fourth-moment sequences}

The next result is important because a pointwise inverse-moment estimate is not enough for a quantitative limit theorem.

\begin{lemma}\label{lem:Laplace-Hermite}
	Let $Z=H_{n-1}(\xi)^2$.  For every $\alpha\in(0,\dfrac{1}{2})$, there exists $c=c(n,\alpha)>0$ such that
	\begin{equation}\label{eq:Laplace-Hermite}
		\mathbb{E}\left[e^{-tZ}\right]
		\leq
		(1+ct)^{-\alpha},
		\qquad t\geq0.
	\end{equation}
\end{lemma}

\begin{proof}
	Lemma~\ref{lem:Hermite-small-ball} gives
	\begin{equation*}
		\mathbb{P}(Z\leq u)\leq C_nu^{\frac{1}{2}},
		\qquad 0<u\leq1.
	\end{equation*}
	Integration by parts yields
	\begin{equation*}
		\mathbb{E}\left[e^{-tZ}\right]
		=
		t\int_0^\infty e^{-tu}\mathbb{P}(Z\leq u)\,du,
	\end{equation*}
	because $\mathbb{P}(Z=0)=0$.  Hence the left-hand side is $O(t^{-\frac{1}{2}})$ as $t\to\infty$.  Near zero,
	\begin{equation*}
		\mathbb{E}\left[e^{-tZ}\right]
		=
		1-t\mathbb{E}[Z]+o(t),
		\qquad
		\mathbb{E}[Z]=(n-1)!>0.
	\end{equation*}
	Choose $c>0$ small enough that the desired inequality holds near zero.  Since $\alpha<\dfrac{1}{2}$, it also holds for all sufficiently large $t$.  Decreasing $c$ once more handles the remaining compact interval.
\end{proof}

\subsection{Proof of Theorem~\ref{thm:uniform-inverse}}

Write the sequence in the form
\begin{equation}\label{eq:Fm-infinite-sum}
	F_m
	=
	\sum_{k\geq1}\lambda_{m,k}H_n(\xi_{m,k}),
	\qquad
	n!\sum_{k\geq1}\lambda_{m,k}^2=1.
\end{equation}
Here $\xi_{m,k}=X(e_{m,k})$, where $(e_{m,k})_{k\geq1}$ is orthonormal for each fixed $m$; no common diagonalizing basis is assumed as $m$ varies.  The sum in \eqref{eq:Fm-infinite-sum} is over all positive integers $k$, whereas $m$ labels the approximating sequence.  When a finite truncation is needed, we write
\begin{equation*}
	F_m^{(M)}=\sum_{k=1}^{M}\lambda_{m,k}H_n(\xi_{m,k}).
\end{equation*}

\begin{proof}
	Write
	\begin{equation*}
		\gamma_{F_m}
		=
		\sum_{k\geq1}a_{m,k}Z_{m,k},
		\qquad
		a_{m,k}=n^2\lambda_{m,k}^2,
		\qquad
		Z_{m,k}=H_{n-1}(\xi_{m,k})^2.
	\end{equation*}
	Set
	\begin{equation*}
		S_n:=\sum_{k\geq1}a_{m,k}=\frac{n^2}{n!}=\frac{n}{(n-1)!},
		\qquad
		b_m:=\sup_{k\geq1}a_{m,k}.
	\end{equation*}
	By Proposition~\ref{prop:diagonal-fourth} and the normalization,
	\begin{equation*}
		b_m
		=n^2\sup_k\lambda_{m,k}^2
		\longrightarrow0.
	\end{equation*}
	Indeed,
	\begin{equation*}
		\sup_k\left\lvert \lambda_{m,k}\right\rvert^4
		\leq
		\sum_k\lambda_{m,k}^4
		=
		\frac{\kappa_4(F_m)}{c_n}.
	\end{equation*}

	Fix $\alpha\in(0,\dfrac{1}{2})$.  Independence and Lemma~\ref{lem:Laplace-Hermite} give
	\begin{equation*}
		\mathbb{E}\left[e^{-t\gamma_{F_m}}\right]
		\leq
		\prod_{k\geq1}\left(1+cta_{m,k}\right)^{-\alpha}.
	\end{equation*}
	The function $x\mapsto\log(1+ctx)$ is concave and vanishes at zero.  Thus, for $0\leq x\leq b_m$,
	\begin{equation*}
		\log(1+ctx)
		\geq
		\frac{x}{b_m}\log(1+ctb_m).
	\end{equation*}
	Consequently,
	\begin{equation}\label{eq:Laplace-product}
		\mathbb{E}\left[e^{-t\gamma_{F_m}}\right]
		\leq
		\left(1+ctb_m\right)^{-\frac{\alpha S_n}{b_m}}.
	\end{equation}

	Fix $q>0$.  The Mellin--Laplace identity gives
	\begin{equation*}
		\mathbb{E}\left[\gamma_{F_m}^{-q}\right]
		=
		\frac{1}{\Gamma(q)}
		\int_0^\infty t^{q-1}
		\mathbb{E}\left[e^{-t\gamma_{F_m}}\right]dt.
	\end{equation*}
	Let
	\begin{equation*}
		A_m:=\frac{\alpha S_n}{b_m}.
	\end{equation*}
	Because $b_m\to0$ and $S_n>0$, one has $A_m\to\infty$.
	For all sufficiently large $m$, $A_m>q$.  Formula \eqref{eq:Laplace-product} and the beta integral imply
	\begin{eqnarray}
		\mathbb{E}\left[\gamma_{F_m}^{-q}\right]
		&\leq&
		\frac{1}{\Gamma(q)}
		\int_0^\infty t^{q-1}(1+ctb_m)^{-A_m}\,dt
		\label{eq:Gamma-ratio-bound}\\
		&=&
		(cb_m)^{-q}\frac{\Gamma(A_m-q)}{\Gamma(A_m)}.
	\nonumber
	\end{eqnarray}
	Stirling's formula gives
	\begin{eqnarray*}
		\frac{\Gamma(A_m-q)}{\Gamma(A_m)}
		&=&
		\frac{(A_m-q)^{A_m-q-\frac{1}{2}}e^{-(A_m-q)}}
		{A_m^{A_m-\frac{1}{2}}e^{-A_m}}
		\left(1+O\left(A_m^{-1}\right)\right)\\
		&=&
		A_m^{-q}e^q
		\left(1-\frac{q}{A_m}\right)^{A_m-q-\frac{1}{2}}
		\left(1+O\left(A_m^{-1}\right)\right)\\
		&=&
		A_m^{-q}\left(1+O\left(A_m^{-1}\right)\right),
	\end{eqnarray*}
	Since $A_m=\dfrac{\alpha S_n}{b_m}$,
	\begin{equation*}
		(cb_m)^{-q}A_m^{-q}
		=(c\alpha S_n)^{-q}.
	\end{equation*}
	This proves \eqref{eq:uniform-negative}.

	It remains to prove \eqref{eq:inverse-limit}.  By Hermite orthogonality and the normalization,
	\begin{eqnarray*}
		\mathbb{E}[\gamma_{F_m}]
		&=&
		n^2(n-1)!\sum_k\lambda_{m,k}^2
	=n,\\
		{\rm Var}(\gamma_{F_m})
		&=&
		n^4{\rm Var}\left(H_{n-1}(\xi)^2\right)
		\sum_k\lambda_{m,k}^4
	\longrightarrow0.
	\end{eqnarray*}
	Thus $\gamma_{F_m}\to n$ in $L^2$, hence in probability.  The continuous mapping theorem yields
	\begin{equation*}
		\gamma_{F_m}^{-q}\longrightarrow n^{-q}
		\qquad\text{in probability}.
	\end{equation*}
	Choose $\delta>0$.  Applying \eqref{eq:uniform-negative} with $q+\delta$ gives
	\begin{equation*}
		\sup_{m\geq m_1}\mathbb{E}\left[\gamma_{F_m}^{-(q+\delta)}\right]<\infty.
	\end{equation*}
	Hence $X_m:=\gamma_{F_m}^{-q}$ is bounded in
	\begin{equation*}
		L^{1+\frac{\delta}{q}},
	\end{equation*}
	on the tail $m\geq m_1$, and is therefore uniformly integrable there.  More explicitly, with $r=1+\dfrac{\delta}{q}>1$,
	\begin{equation*}
		\sup_{m\geq m_1}\mathbb{E}\left[X_m\mathbf{1}_{\{X_m>R\}}\right]
		\leq
		R^{1-r}\sup_{m\geq m_1}\mathbb{E}[X_m^r]
		\longrightarrow0
	\end{equation*}
	as $R\to\infty$.  Vitali's convergence theorem now gives
	\begin{equation*}
		\mathbb{E}\left\lvert \gamma_{F_m}^{-q}-n^{-q}\right\rvert\longrightarrow0,
	\end{equation*}
	which proves \eqref{eq:inverse-limit}.
\end{proof}

\subsection{Density estimates from the coarea formula}

Assume that $F$ is smoothly nondegenerate.  Then $\delta_y(F)$ is a Watanabe distribution and the coarea formula gives the following representation of the density of $F$ \cite{LI04,N06}:
\begin{equation}\label{eq:coarea-density}
	p_F(y)
	=
	\mathbb{E}\left[\delta_y(F)\right]
	=
	\int_{F^{-1}(y)}\gamma_F^{-\frac{1}{2}}\,da_y,
\end{equation}
where $da_y$ denotes the corresponding Malliavin surface measure on the level set $F^{-1}(y)$.

For a smooth compactly supported function $\varphi$, Malliavin integration by parts gives
\begin{equation}\label{eq:density-ibp}
	\mathbb{E}\left[\varphi'(F)\right]
	=
	\mathbb{E}\left[\varphi(F)H_F\right],
	\qquad
	H_F:=\delta\left(\frac{DF}{\gamma_F}\right).
\end{equation}
Since $-LF=nF$, the product rule for the divergence operator yields
\begin{equation}\label{eq:HF-formula}
	H_F
	=
	\frac{nF}{\gamma_F}
	+
	\frac{\left\langle D\gamma_F,DF\right\rangle_{\mathbb{H}}}{\gamma_F^2}.
\end{equation}
Moreover,
\begin{equation*}
	D\gamma_F
	=
	2\left\langle D^2F,DF\right\rangle_{\mathbb{H}},
\end{equation*}
where contraction over one Hilbert-space coordinate produces an element of $\mathbb{H}$.  Therefore,
\begin{equation*}
	\left\lvert \left\langle D\gamma_F,DF\right\rangle_{\mathbb{H}}\right\rvert
	\leq
	2\left\lVert D^2F\right\rVert_{\mathbb{H}^{\otimes2}}\left\lVert DF\right\rVert_{\mathbb{H}}^2
	=
	2\left\lVert D^2F\right\rVert_{\mathbb{H}^{\otimes2}}\gamma_F,
\end{equation*}
and hence
\begin{equation}\label{eq:HF-bound}
	\left\lvert H_F\right\rvert
	\leq
	\frac{n\left\lvert F\right\rvert+2\left\lVert D^2F\right\rVert_{\mathbb{H}^{\otimes2}}}{\gamma_F}.
\end{equation}

Comparing the two sides of \eqref{eq:density-ibp} through the density of $F$ gives, in the distributional sense,
\begin{equation}\label{eq:p-prime-conditional}
	-\nabla \log p_F(x)
	=
	\mathbb{E}\left[H_F\mid F=x\right].
\end{equation}
Consequently, conditional Jensen gives
\begin{eqnarray*}
	\left\lVert p_F'\right\rVert_{L^1(\mathbb{R})}
	&=&
	\int_{\mathbb{R}}
	\left\lvert \mathbb{E}[H_F\mid F=x]\right\rvert p_F(x)\,dx\\
	&\leq&
	\mathbb{E}\left\lvert H_F\right\rvert.
\end{eqnarray*}
	Let $q>1$ and let $r=\dfrac{q}{q-1}$.  By \eqref{eq:HF-bound} and H\"older's inequality,
	\begin{eqnarray*}
		\left\lVert p_F'\right\rVert_{L^1(\mathbb{R})}
		&\leq&
		n\left\lVert F\right\rVert_{L^r}\left\lVert \gamma_F^{-1}\right\rVert_{L^q}\\
		&&\quad+
		2\left\lVert \left\lVert D^2F\right\rVert_{\mathbb{H}^{\otimes2}}\right\rVert_{L^r}
		\left\lVert \gamma_F^{-1}\right\rVert_{L^q}.
	\end{eqnarray*}
	Because $F$ belongs to the fixed $n$th chaos and has unit variance, hypercontractivity gives
	\begin{equation*}
		\left\lVert F\right\rVert_{L^r}
		+
		\left\lVert \left\lVert D^2F\right\rVert_{\mathbb{H}^{\otimes2}}\right\rVert_{L^r}
		\leq
		C_{n,q}.
\end{equation*}
Here we used the identity
\begin{equation*}
	\mathbb{E}\left\lVert D^2F\right\rVert_{\mathbb{H}^{\otimes2}}^2
	=
	n(n-1)\mathbb{E}[F^2]
	=
	n(n-1).
\end{equation*}
It follows that
\begin{equation}\label{eq:p-prime-L1}
	\left\lVert p_F'\right\rVert_{L^1(\mathbb{R})}
	\leq
		C_{n,q}\left(1+\left\lVert \gamma_F^{-1}\right\rVert_{L^q}\right).
\end{equation}

\subsection{A weighted density estimate}

Denote by $\mu_F$ the distribution of $F$ and by $\gamma$ the standard Gaussian distribution.  Let $\phi$ and $\Phi$ be the density and distribution function of $\gamma$, respectively.  Put
\begin{equation*}
	g(x):=p_F(x)-\phi(x),
	\qquad
	G(x):=\int_{-\infty}^{x}g(u)\,du.
\end{equation*}
Then
\begin{equation*}
	\left\lVert G\right\rVert_\infty=d_{\mathrm K}\left(\mu_F,\gamma\right).
\end{equation*}
Define
\begin{equation*}
	B_F:=1+\left\lVert p_F'\right\rVert_{L^1(\mathbb{R})}.
\end{equation*}
Since $\left\lVert \phi'\right\rVert_{L^1}<\infty$,
\begin{equation}\label{eq:g-prime}
	\left\lVert g'\right\rVert_{L^1(\mathbb{R})}
	\leq
	C B_F.
\end{equation}

\begin{lemma}\label{lem:weighted-L1}
	For every $p\geq1$ and $\eta>0$, there exists $C_{p,\eta}>0$ such that
	\begin{equation}\label{eq:weighted-L1}
		\int_{\mathbb{R}}\left\lvert x\right\rvert^p\left\lvert p_F(x)-\phi(x)\right\rvert\,dx
		\leq
		C_{p,\eta}
		\left(d_{\mathrm K}\left(\mu_F,\gamma\right)B_F\right)^{
		\frac{\eta}{2(p+\eta+\frac{1}{2})}}.
	\end{equation}
\end{lemma}

\begin{proof}
	Fix an interval $I=[u,v]$ of length at most $h$, and let
	\begin{equation*}
		g_I:=\frac{1}{\left\lvert I\right\rvert}\int_Ig(x)\,dx.
	\end{equation*}
	For every $x\in I$,
	\begin{equation*}
		\left\lvert g(x)-g_I\right\rvert
		\leq
		\frac{1}{\left\lvert I\right\rvert}\int_I\left\lvert g(x)-g(y)\right\rvert\,dy
		\leq
		\operatorname{Var}_I(g),
	\end{equation*}
	where $\operatorname{Var}_I(g)$ is the total variation of $g$ on $I$.  Hence
	\begin{equation}\label{eq:interval-variation}
		\int_I\left\lvert g(x)\right\rvert\,dx
		\leq
		\left\lvert \int_Ig(x)\,dx\right\rvert
		+
	h\operatorname{Var}_I(g).
	\end{equation}
	Because $g$ is absolutely continuous, every partition $u=x_0<\cdots<x_J=v$ satisfies
	\begin{equation*}
		\sum_{j=1}^{J}\left\lvert g(x_j)-g(x_{j-1})\right\rvert
		\leq
		\int_I\left\lvert g'(x)\right\rvert\,dx.
	\end{equation*}
	Taking the supremum over partitions yields
	\begin{equation}\label{eq:variation-derivative}
		\operatorname{Var}_I(g)
		\leq
		\int_I\left\lvert g'(x)\right\rvert\,dx.
	\end{equation}
	Also,
	\begin{equation*}
		\left\lvert \int_Ig(x)\,dx\right\rvert
		=
		\left\lvert G(v)-G(u)\right\rvert
		\leq
		2d_{\mathrm K}\left(\mu_F,\gamma\right).
	\end{equation*}

	Fix $R\geq1$ and $0<h\leq R$. Partition $[-R,R]$ into at most $C R/h$ intervals, each of length at most $h$. Summing \eqref{eq:interval-variation} and using \eqref{eq:variation-derivative}--\eqref{eq:g-prime} gives
	\begin{equation*}
		\int_{-R}^{R}\left\lvert g(x)\right\rvert\,dx
		\leq
		C\left(
		\frac{R}{h}d_{\mathrm K}\left(\mu_F,\gamma\right)
		+
		hB_F
		\right).
	\end{equation*}
	When $d_{\mathrm K}\left(\mu_F,\gamma\right)>0$, choose
	\begin{equation*}
		h=\left(\frac{R\,d_{\mathrm K}\left(\mu_F,\gamma\right)}{B_F}\right)^{\frac{1}{2}}.
	\end{equation*}
	Since $d_{\mathrm K}\left(\mu_F,\gamma\right)\leq1$ and $B_F\geq1$, one has $h\leq R$.  Thus
	\begin{equation}\label{eq:weighted-inside}
		\int_{\left\lvert x\right\rvert\leq R}\left\lvert x\right\rvert^p\left\lvert g(x)\right\rvert\,dx
		\leq
		CR^{p+\frac{1}{2}}
		\left(d_{\mathrm K}\left(\mu_F,\gamma\right)B_F\right)^{\frac{1}{2}}.
	\end{equation}
	If the Kolmogorov distance is zero, then the two laws coincide and the conclusion is immediate.

	For the tail and every $\eta>0$,
	\begin{eqnarray}
		\int_{\left\lvert x\right\rvert>R}\left\lvert x\right\rvert^p\left\lvert g(x)\right\rvert\,dx
		&\leq&
		\mathbb{E}\left[\left\lvert F\right\rvert^p\mathbf{1}_{\{\left\lvert F\right\rvert>R\}}\right]
		+
		\mathbb{E}\left[\left\lvert N\right\rvert^p\mathbf{1}_{\{\left\lvert N\right\rvert>R\}}\right]
		\label{eq:weighted-outside}\\
		&\leq&
		R^{-\eta}\left(\mathbb{E}\left\lvert F\right\rvert^{p+\eta}+\mathbb{E}\left\lvert N\right\rvert^{p+\eta}\right)
		\nonumber\\
		&\leq&
		C_{p,\eta}R^{-\eta}.
	\nonumber
	\end{eqnarray}
	The last step uses hypercontractivity and $\mathbb{E}[F^2]=1$.  Write
	\begin{equation*}
		A:=d_{\mathrm K}\left(\mu_F,\gamma\right)B_F.
	\end{equation*}
	If $0<A<1$, balance \eqref{eq:weighted-inside} and \eqref{eq:weighted-outside} by taking
	\begin{equation*}
		R=A^{-\frac{1}{2(p+\eta+\frac{1}{2})}}.
	\end{equation*}
	This gives \eqref{eq:weighted-L1}.  If $A\geq1$, the result follows after increasing the constant, since the weighted integral is bounded by $\mathbb{E}\left\lvert F\right\rvert^p+\mathbb{E}\left\lvert N\right\rvert^p$.
\end{proof}

\subsection{Proof of Theorem~\ref{thm:diagonal-non-Stein} and Corollary~\ref{cor:diagonal-convergence}}

Before proving Theorem ~\ref{thm:diagonal-non-Stein}, we introduce the following useful lemma.
\begin{lemma}\cite{VIL09}\label{lem:weighted-TV-Wp}
	Let $\mu(dx)=u(x)\,dx$ and $\nu(dx)=v(x)\,dx$ be probability measures on $\mathbb{R}$ with finite $p$-th moments.  Then
	\begin{equation}\label{eq:weighted-TV-Wp}
		W_p^p(\mu,\nu)
		\leq
		2^{p-1}\int_{\mathbb{R}}\left\lvert x\right\rvert^p\left\lvert u(x)-v(x)\right\rvert\,dx.
	\end{equation}
\end{lemma}
 
 Now we give the proof of Theorem~\ref{thm:diagonal-non-Stein}.
\begin{proof}
	Apply Lemma~\ref{lem:weighted-TV-Wp} with $u=p_F$ and $v=\phi$.  Since $\mu_F(dx)=p_F(x)\,dx$ and $\gamma(dx)=\phi(x)\,dx$, we first obtain
	\begin{equation*}
		W_p^p\left(\mu_F,\gamma\right)
		\leq
		2^{p-1}
		\int_{\mathbb{R}}\left\lvert x\right\rvert^p\left\lvert p_F(x)-\phi(x)\right\rvert\,dx.
	\end{equation*}
	Lemma~\ref{lem:weighted-L1} bounds the integral on the right-hand side.  Taking the $p$th root gives
	\begin{eqnarray*}
		W_p\left(\mu_F,\gamma\right)
		&\leq&
		C_{p,\eta,n}
		\left(d_{\mathrm K}\left(\mu_F,\gamma\right)B_F\right)^{
		\frac{\eta}{2p\left(p+\eta+\frac{1}{2}\right)}}.
	\end{eqnarray*}
	By Proposition~\ref{prop:diagonal-BE} and \eqref{eq:p-prime-L1},
	\begin{eqnarray*}
		d_{\mathrm K}\left(\mu_F,\gamma\right)
		&\leq&
		C_n\sqrt{\kappa_4(F)},
		\\
		B_F
		&\leq&
		C_{n,q}
		\left(
		1+
		\left(\mathbb{E}\left[\gamma_F^{-q}\right]\right)^{\frac{1}{q}}
		\right).
	\end{eqnarray*}
	Substituting these two estimates into the preceding display gives
	\begin{equation*}
		W_p\left(\mu_F,\gamma\right)
		\leq
		C_{p,\eta,n,q}
		\left(
		1+
		\left(\mathbb{E}\left[\gamma_F^{-q}\right]\right)^{\frac{1}{q}}
		\right)^{
		\frac{\eta}{2p\left(p+\eta+\frac{1}{2}\right)}}
		\left(\sqrt{\kappa_4(F)}\right)^{
		\frac{\eta}{2p\left(p+\eta+\frac{1}{2}\right)}}.
	\end{equation*}
	Finally,
	\begin{equation*}
		\left(\sqrt{\kappa_4(F)}\right)^{
		\frac{\eta}{2p\left(p+\eta+\frac{1}{2}\right)}}
		=
		\kappa_4(F)^{
		\frac{\eta}{4p\left(p+\eta+\frac{1}{2}\right)}}.
	\end{equation*}
	This is exactly \eqref{eq:diagonal-main}.
\end{proof}

Consequensely, we deduce the proof of Corollary~\ref{cor:diagonal-convergence}.

\begin{proof}
	Fix any $q>1$.  Theorem~\ref{thm:uniform-inverse} gives
	\begin{equation*}
		\sup_{m\geq m_0}\mathbb{E}\left[\gamma_{F_m}^{-q}\right]<\infty.
	\end{equation*}
	Consequently, there exists a finite constant $K_q$ such that
	\begin{equation*}
		1+
		\left(\mathbb{E}\left[\gamma_{F_m}^{-q}\right]\right)^{\frac{1}{q}}
		\leq
		K_q,
		\qquad m\geq m_0.
	\end{equation*}
	The infinite-support assumption and Corollary~\ref{cor:smooth-nondegenerate} imply that every $F_m$ is smoothly nondegenerate.  Notice that Theorem~\ref{thm:uniform-inverse} alone would not justify this assertion: its index $m_0$ is allowed to depend on the prescribed inverse-moment order, whereas smooth nondegeneracy requires all negative moments for each fixed $m$.

	We may therefore apply Theorem~\ref{thm:diagonal-non-Stein} to every $m\geq m_0$.  The factor involving the negative moment is bounded by a constant depending only on $q$ and is absorbed into $C_{p,\eta,n}$.  This gives \eqref{eq:diagonal-sequence-rate}.  Since
	\begin{equation*}
		\frac{\eta}{4p\left(p+\eta+\frac{1}{2}\right)}>0,
	\end{equation*}
	the assumption $\kappa_4(F_m)\to0$ implies $W_p\left(\mu_m,\gamma\right)\to0$.
\end{proof}

\end{document}